\documentclass[10pt,twocolumn,twoside,journal]{IEEEtran}

\title{Quantitative Sensitivity Bounds for Nonlinear Programming and Time-varying Optimization}
\author{Irina Suboti\'c, \IEEEmembership{Student Member, IEEE}, Adrian Hauswirth, \IEEEmembership{Student Member, IEEE}, and Florian D\"orfler, \IEEEmembership{Member, IEEE}
\thanks{This work was supported in part by ETH Z\"urich funds and in part by the Swiss Federal Office of Energy grant \#SI/501707 GREAT.}
\thanks{The authors are with the Department of Information Technology and Eletrical Engineering, ETH Z\"urich, Physikstrasse 3, 8092 Z\"urich, Switzerland.}
\thanks{Email: \{subotici,hadrian,dorfler\}@ethz.ch.}}

\usepackage{cite}

\usepackage{graphicx}

\usepackage{amsmath}
\usepackage{amssymb}
\usepackage{amsfonts}
\usepackage{enumerate}
\usepackage{booktabs}
\usepackage{enumitem}
\usepackage[caption=false,font=footnotesize]{subfig}

\usepackage{amsthm}
\usepackage{thmtools}

\newtheorem{theorem}{Theorem}
\newtheorem{lemma}{Lemma}
\newtheorem{definition}{Definition}
\newtheorem{assumption}{Assumption}
\newtheorem{proposition}[theorem]{Proposition}
\newtheorem{corollary}[theorem]{Corollary}
\theoremstyle{remark}

\newcommand\oprocendsymbol{\hbox{\small $\blacksquare$}}
\newcommand\oprocend{\relax\ifmmode\else\unskip\hfill\fi\oprocendsymbol}

\DeclareMathAlphabet{\mymathbb}{U}{bbold}{m}{n}

\newcommand{\bbone}{\mymathbb{1}}
\newcommand{\bbzero}{\mymathbb{0}}

\DeclareMathOperator{\diag}{diag}
\DeclareMathOperator*{\minimize}{minimize}
\DeclareMathOperator*{\st}{subject~to}
\DeclareMathOperator{\cl}{cl}

\newcommand{\R}{\mathbb{R}}
\newcommand{\N}{\mathbb{N}}
\newcommand{\calX}{{\mathcal{X}}}
\newcommand{\calN}{\mathcal{N}}
\newcommand{\bbI}{\mathbb{I}}
\newcommand{\bbB}{\mathbb{B}}

\newcommand{\bfR}{\mathbf{R}}
\newcommand{\bfI}{\mathbf{I}}

\newcommand{\bfIS}{\mathbf{I}^{\mathbf{s}}}

\newcommand{\T}{\top}
\newcommand{\coloneqq}{:=}
\providecommand{\norm}[1]{\lVert#1\rVert}
\providecommand{\abs}[1]{\lvert#1\rvert}

\newcommand{\txtold}{\tilde{\xi}(t)}
\newcommand{\txt}{t}
\newcommand{\txs}{\tilde{x}^\star}

\newenvironment{smallbmatrix}
  {\left[\begin{smallmatrix}}
  {\end{smallmatrix}\right]}

\newcommand{\grad}{\boldsymbol \nabla}

\usepackage{cite}
\usepackage[capitalize,nameinlink,compress]{cleveref}
\crefname{appendix}{Appendix}{Appendices} 
\crefname{figure}{Figure}{Figures} 
\crefname{line}{line}{lines} 
\crefname{claim}{Claim}{Claims} 
\crefname{equation}{}{} 
\crefname{problem}{Problem}{Problems}
\crefname{proposition}{Proposition}{Propositions}
\crefname{assumption}{Assumption}{Assumptions}
\crefname{remarkx}{Remark}{Remarks}

\newlist{romcases}{enumerate}{1}
\setlist[romcases,1]{ref=\roman*, label=(\roman*)}
\creflabelformat{romcasesi}{{\it(#1)}}
\crefname{romcasesi}{}{}

\begin{document}

\maketitle

\begin{abstract}
    Inspired by classical sensitivity results for nonlinear optimization, we derive and discuss new quantitative bounds to characterize the solution map and dual variables of a parametrized nonlinear program. In particular, we derive explicit expressions for the local and global Lipschitz constants of the solution map of  non-convex or convex optimization problems, respectively.
    Our results are geared towards the study of time-varying optimization problems which are commonplace in various applications of online optimization, including power systems, robotics, signal processing and more. In this context, our results can be used to bound the rate of change of the optimizer.
    To illustrate the use of our sensitivity bounds we generalize existing arguments to quantify the tracking performance of continuous-time, monotone running algorithms. Further, we introduce a new continuous-time running algorithm for time-varying constrained optimization which we model as a so-called perturbed sweeping process. For this discontinuous scheme we establish an explicit bound on the asymptotic solution tracking for a class of convex problems.
\end{abstract}


\section{Introduction}

\IEEEPARstart{S}{olving} optimization problems that vary over time is a widely encountered challenge referred to as \emph{time-varying} or \emph{non-stationary optimization}\cite{simonettoTimeVaryingConvexOptimization2017,tangDistributedAlgorithmTimevarying2017}, or as online optimization in \emph{dynamic environments}\cite{lesage-landryPredictiveOnlineConvex2020, mokhtariOnlineOptimizationDynamic2016}.

Notable applications where time-varying optimization problems arise include feedback-based optimization~\cite{hauswirthAntiWindupApproximationsOblique2020,tangRunningPrimalDualGradient2018,hauswirthTimevaryingProjectedDynamical2018,colombinoOnlineOptimizationFeedback2019,bernsteinOnlineOptimizationFeedback2018} (e.g. for power systems \cite{liuCoordinateDescentAlgorithmTracking2018,dallaneseOptimalPowerFlow2018,hauswirthOnlineOptimizationClosed2017}), real-time iterations \cite{zavalaRealTimeNonlinearOptimization2010,diehlRealtimeOptimizationNonlinear2002,fellerStabilizingIterationScheme2017} in model predictive control, formation and collision avoidance in robotics~\cite{rahiliDistributedConvexOptimization2015,sunDistributedTimeVaryingFormation2018}, and many more.
Roughly speaking, for these problems we can distinguish between \emph{running} algorithms\cite{popkovGradientMethodsNonstationary2005,mokhtariOnlineOptimizationDynamic2016}, that do not anticipate the change in an optimizer, as opposed to \emph{predictive} schemes that exploit some information about the change of the optimization problem \cite{simonettoPredictionCorrectionAlgorithmsTimeVarying2017, lesage-landryPredictiveOnlineConvex2020, fazlyabPredictionCorrectionInteriorPointMethod2018, simonettoDecentralizedPredictionCorrectionMethods2016}.

Depending the application, the optimization problem might be unconstrained with a time-varying cost function~\cite{popkovGradientMethodsNonstationary2005,maddenBoundsTrackingError2020}, constrained to a stationary set but with time-varying cost function \cite{mokhtariOnlineOptimizationDynamic2016,zinkevichOnlineConvexProgramming2003}, or have time-varying constraints as well\cite{zhangOnlineProximaladmmTimevarying2020}. If constraints are time-varying it is possible to eliminate equality constraints or apply barrier functions \cite{fazlyabInteriorPointMethod2016,fazlyabPredictionCorrectionInteriorPointMethod2018} to the objective to arrive at an unconstrained problem.

All of these setups require the solution of the time-varying problem to be a well-behaved function of time. In other words, an optimizer should be isolated, if not unique, and vary slowly so that an online algorithm is able to approximately track it. We refer to the asymptotic distance between the true optimizer and the algorithm state (or iterate) as \emph{tracking performance}. To analyze the tracking performance, most works have focused on convex problems, with strongly convex objectives that are guaranteed to admit a unique global optimizer. Moreover, many works establish tracking performance under the assumption that the rate of change of the optimizer is bounded. However, apart from special cases, it is unclear if and how this assumption is satisfied, and which rates are achievable.

Therefore, in this paper, we study how the rate of change of optimizers can be bounded for a general class of nonlinear optimization problems. More precisely, using the classical sensitivity results for parametrized, nonlinear optimization problems
\cite{fiaccoSensitivityAnalysisNonlinear1976, fiaccoSensitivityStabilityAnalysis1990, robinsonPerturbedKuhnTuckerPoints1974, jittorntrumSolutionPointDifferentiability1984} we establish sufficient conditions under which the instantaneous optimizer and its dual variables are continuous functions of a perturbation (not necessarily time), with a quantifiable sensitivity. Moreover, we provide and discuss special cases of constrained parametrized optimization problems for which the sensitivities of optimizer and dual variables take a simple and easy-to-interpret form.

As a second contribution, in view of the renewed interest in continuous-time optimization \cite{helmkeOptimizationDynamicalSystems1996,goebelHybridDynamicalSystems2012,muehlebachDynamicSystemPerspective2019}, we show the usefulness of our sensitivity results by quantifying the tracking performance of continuous-time running algorithms for time-varying optimization. In particular, we generalize and formalize a monotonicity-like property that can be used to provide asymptotic bounds on the tracking error of, almost arbitrary, well-defined continuous-time schemes.

Finally, as a third contribution, we illustrate our results by introducing a novel continuous-time algorithm, termed \emph{sweeping gradient flow} and modeled as a {perturbed sweeping process}  \cite{edmondRelaxationOptimalControl2005,thibaultSweepingProcessRegular2003}, to track the solution of a convex optimization problem with time-varying objective and constraints.
The proposed scheme is inherently discontinuous and defined by a differential inclusion. However, the scheme is of conceptual importance since it can be interpreted as the continuous-time limit of an online projected gradient descent \cite{mokhtariOnlineOptimizationDynamic2016,zinkevichOnlineConvexProgramming2003}. Moreover, sweeping processes of this type, which essentially model a continuous projection onto a time-varying set, have been considered in feedback optimization~\cite{tangRunningPrimalDualGradient2018,hauswirthTimevaryingProjectedDynamical2018,hauswirthAntiWindupApproximationsOblique2020} where constraint enforcement may rely on input saturation in physical plants.

The paper is structured as follows: In \cref{sec:prelim} we fix the notation and recall standard optimality conditions for nonlinear programming, as well as results from sensitivity analysis. \Cref{sec:lipschitz.conitnulity.of.solution.maps} deals with quantifying the local differentiability and Lipschitz continuity of the solution maps for nonlinear programs. Then, in \cref{sec:assumptions.for.global.solution.maps} we introduce and discuss assumptions that allow for global statements about solution maps and their sensitivities. These findings are discussed in \cref{sec:disc}, using illustrative special cases. \Cref{sec:tracking.performance} provides a general result to quantify the tracking performance of continuous-time algorithms under a monotonicity-like property. Next, \cref{sec:sweeping.gradient.descent} introduces a ``sweeping gradient flow'' as a continuous-time scheme for the time-varying optimization and discusses its tracking guarantees for a special class of convex problems. Finally, \cref{sec:numerical.examples} provides illustrative examples validating our theoretical results. Finally, \cref{sec:conc} summarizes our findings and discusses open problems.

\section{Preliminaries}\label{sec:prelim}

 We consider $\R^n$ with the usual inner product $\left\langle \cdot, \cdot \right\rangle$ and 2-norm $\| \cdot \|$.
 The non-negative orthant is denoted by $\R^n_{\geq 0}$.
 The identity and zero matrices of appropriate dimension are written as $\bbI$ and $\bbzero$, respectively.
 For a vector $v \in \R^n$ and an index set $\bfI \subset \{1, \ldots, n \}$, $v_\bfI \in \R^{| \bfI |}$ denotes the vector corresponding to the $\bfI$-th components of $v$. If $A \in \R^{n \times m}$ is a matrix, $A_\bfI$ denotes the submatrix consisting of the $\bfI$-th rows of $A$.
 Furthermore, $A^\mathsf{T}$ denotes its transpose of $A$ and $\|A\| := \sup_{x \in \R^m} \| A x \| / \| x\|$ its operator norm.
 If $A \in \R^{n \times n}$ is symmetric, $\lambda^{\max}_A$ and $\lambda^{\min}_A$ denote the maximum and minimum eigenvalues of $A$, respectively. We write $A\succ 0$  $(A \succeq 0)$ if $A$ is positive (semi) definite.
 For a matrix $B \in \R^{n \times m}$, $\sigma^{\max}_B$ and $\sigma^{\min}_B$ are the maximum and minimum singular values of $B$, respectively. Namely, we have $\sigma^{\max}_B = \| B \|$.

  A map $F: \calX \rightarrow \R^m$ with $\calX \subset \R^n$ is \emph{$\ell$-Lipschitz} if
 \begin{align}\label{eq:lip_cond}
    \| F(y) - F(x) \| \leq \ell\| y - x \|
 \end{align}
 holds for all $y, x \in \calX$, and $F$ is \emph{locally $\ell'$-Lipschitz at $x^\prime \in \calX$} if for every $\ell \geq  \ell'$ there exists a (relative) neighborhood $\calN$ of $x^\prime$ such that \eqref{eq:lip_cond} holds for all $x, y \in \calN$. In other words, $\ell'$ is the largest lower bound on $\ell$ such that \eqref{eq:lip_cond} is satisfied, i.e., $\ell' = \sup_{y' \rightarrow x} \| F(y) - F(x) \| / \| y - x \|$.

 Derivatives are understood in the sense of Fr\'echet, and we denote the Jacobian matrix of $F$ at $x$ by $\nabla F(x) \in \R^{m \times n}$.
 If $F$ is differentiable at $x$ then it is locally $\ell$-Lipschitz at $x$ with $\ell = \| \nabla F(x) \|$. The Jacobian at $x$, with respect to a subset of variables $x_1 \in \R^{n_1}$ is denoted by $\nabla_{x_1} F(x) \in \R^{m\times n_1}$.

 If $m = 1$, we call the column vector $\grad F(x) := \nabla F(x)^\T$ the \emph{gradient} of $F$, and $\nabla_{xx}^2 \in \R^{n \times n}$ denotes the Hessian of $F$ at $x$. However, when differentiating $F$ twice, with respect to different subsets of variables $x_1 \in \R^{n_1}$ and $x_2 \in \R^{n_2}$, we define $\nabla^2_{x_1 x_2} F(x) :=\nabla_{x_1}(\grad_{x_2}F)(x) \in \R^{n_2\times n_1}$, in particular, we have $ \nabla^2_{x_2 x_1} F(x)= (\nabla^2_{x_1 x_2} F(x))^\T$.

\subsection{Sensitivity Analysis for Nonlinear Optimization}

Throughout the paper we consider the following parametrized nonlinear optimization problem and special cases thereof:
\begin{align}\label{eq:basic_opt_prob}
    \begin{split}
        \underset{x}{\minimize} \quad & f(x, \xi)\\
        \st \quad & h(x, \xi) = 0\\
                  & g(x, \xi) \leq 0 \, ,
      \end{split}
\end{align}
where $x \in \R^{n}$ is the decision variable, $\xi \in \Xi \subset \R^r$ is a parameter. Further, the objective $f: \R^n \times \Xi \rightarrow \R$ and the constraint functions $h: \R^n \times \Xi \rightarrow \R^p$ and $g: \R^n \times \Xi \rightarrow \R^m$ are twice continuously differentiable in $(x, \xi)$. Further, we define the parameterized feasible set $\calX(\xi) := \{ x \in \R^n \, | \, h(x, \xi) = 0, \, g(x, \xi) \leq 0 \}$.
Given $x \in \calX(\xi)$, let $\bfI_x^\xi := \{ i \, | \, g_i (x, \xi) = 0 \}$ denote the \emph{set of active inequality constraints at $x$ for $\xi \in \Xi$} and $\Bar{\bfI}_x^\xi := \{ i \, | \, g_i (x, \xi ) < 0 \}$ the \emph{set of inactive inequality constraints at $x$ for $\xi \in \Xi$}.

The following standard definitions and results closely follow \cite{jittorntrumSolutionPointDifferentiability1984,jittorntrumSequentialAlgorithmsNonlinear1978,fiaccoSensitivityAnalysisNonlinear1976}, but they can also be found in any good textbook on nonlinear optimization, e.g.,~\cite{bazaraaNonlinearProgrammingTheory2006, luenbergerLinearnonlinearprogramming1984}.

\begin{definition}[LICQ] Given~\eqref{eq:basic_opt_prob}, $\xi \in \Xi$, and $x \in \calX(\xi)$, the \emph{linear independence constraint qualification (LICQ) holds at $x$} if the matrix $\begin{bmatrix} \nabla_x h(x, \xi)^T & \nabla_x g_{\bfI_x^\xi}(x, \xi)^T \end{bmatrix}$ has full rank.
\end{definition}

The \emph{(parametrized) Lagrangian of~\eqref{eq:basic_opt_prob}} is defined as $L(x, \lambda, \mu, \xi) := f(x, \xi) + \lambda^\T h(x, \xi) + \mu^\T g(x, \xi)$ for all $x \in \R^n$, all \emph{Lagrange multipliers} $\mu \in \R_{\geq 0}^m$, and all $\xi \in \Xi$.

\begin{definition}[KKT] Given~\eqref{eq:basic_opt_prob} for some $\xi \in \Xi$, the Karush-Kuhn-Tucker (KKT; first-order optimality) conditions hold at $(x^\star, \lambda^\star, \mu^\star) \in \calX(\xi) \times \R^p \times \R^m_{\geq 0}$ if $\nabla_x L(x^\star, \lambda^\star, \mu^\star, \xi) = 0$ and $\mu_i^\star = 0$ holds for all $i \in \bar{\bfI}^\xi_{x^\star}$.
\end{definition}

Under LICQ (but also under weaker constraint qualifications, such as the Mangasarian-Fromowitz CQ), every optimizer is a KKT point \cite[Th.~4.2.13]{bazaraaNonlinearProgrammingTheory2006}.
However, the uniqueness of $\mu^\star$, holds only under LICQ.
\begin{theorem}
    Given $\xi \in \Xi$, if $x^\star$ is a (local) solution of \eqref{eq:basic_opt_prob} and LICQ holds at $x^\star$, then there exists a unique $(\lambda^\star, \mu^\star) \in \R^p \times \R^{m}_{\geq 0}$ such that the KKT conditions are satisfied at $(x^\star, \lambda^\star, \mu^\star)$.
\end{theorem}

\begin{definition}[SSOSC] Consider \eqref{eq:basic_opt_prob} for $\xi \in \Xi$. The point $(x^\star, \lambda^\star, \mu^\star) \in \calX(\xi) \times \R^p \times \R^m_{\geq 0}$ satisfies the \emph{strong second-order sufficiency condition (SSOSC)} if $y^\T \nabla_{xx}^2 L y  > 0$ holds for all $y \neq 0$ and all $i$ for which $\nabla_x h(x^\star,\xi) y = 0$, $ \nabla_x g_i(x^\star, \xi) y = 0$ with $\mu^\star_i>0$.
\end{definition}

Loosely speaking, the SSOSC guarantees that $\nabla_{xx}L$ is positive definite locally at $(x^\star, \lambda^\star, \mu^\star)$ and along all feasible directions $y$. More precisely, the SSOSC guarantees that a KKT point is a strict optimizer of~\eqref{eq:basic_opt_prob}~\cite[Th.~4.4.2]{bazaraaNonlinearProgrammingTheory2006}:
\begin{theorem}
    Given $\xi \in \Xi$, if $(x^\star, \lambda^\star, \mu^\star)$ satisfies LICQ at $x^\star$, KKT, and SSOSC, then $x^\star$ is a strict local minimizer of~\eqref{eq:basic_opt_prob}.
\end{theorem}

Given $\xi \in \Xi$, we say that $(x^\star, \lambda^\star, \mu^\star) \in \calX(\xi) \times \R^p \times \R^m_{\geq 0}$ is a \emph{regular minimizer of~\eqref{eq:basic_opt_prob}} if LICQ is satisfied at $x^\star$ and $(x^\star, \lambda^\star, \mu^\star)$ satisfies KKT and SSOSC. For regular optimizers, the following sensitivity result holds \cite[Thm. 2.3.3]{jittorntrumSequentialAlgorithmsNonlinear1978}:
\begin{theorem}\label{thm:sens_continuous} Consider~\eqref{eq:basic_opt_prob} and let $f,g,h$ be twice continuously differentiable\footnote{As a slightly less stringent condition for \cref{it:1}, \cref{it:2}, and \cref{it:3} to hold, we require only that $f$ and $g$ are twice continuously differentiable in $x$, and that $f, g, h, \nabla_x f$, $\nabla_x g, \nabla_x h, \nabla^2_{xx} f, \nabla^2_{xx} g$ and $\nabla^2_{xx} h$ are continuous in $\xi$ \cite[Thm.2.3.2]{jittorntrumSequentialAlgorithmsNonlinear1978}.} in $(x, \xi)$. Further, let $\bar{x}^\star$ be a regular minimizer of \eqref{eq:basic_opt_prob} for $\bar{\xi}$ with multiplier $(\bar{\lambda}^\star, \bar{\mu}^\star)$.
Then, on a neighborhood $\calN \subset \Xi$ of $\bar{\xi}$, there exist continuous maps $x^\star: \calN \rightarrow \R^n$, $\lambda^\star: \calN \rightarrow \R^p$, and $\mu^\star: \calN \rightarrow \R_{\geq 0}^m$ such that
\begin{romcases}
    \item \label{it:1} $x^\star(\bar{\xi}) = \bar{x}^\star$, $\lambda^\star(\bar{\xi}) = \bar{\lambda}^\star$, and $\mu^\star(\bar{\xi}) = \bar{\mu}^\star$,
    \item \label{it:2} for all $\xi \in \calN$, LICQ is satisfied for $x^\star(\xi)$, and KKT and SSOSC hold for $(x^\star(\xi), \lambda^\star(\xi), \mu^\star(\xi))$,
    \item \label{it:3} for all $\xi \in \calN$, $x^\star(\xi)$ is a local minimizer for \eqref{eq:basic_opt_prob}, and $(\lambda^\star(\xi), \mu^\star(\xi)$ is the corresponding Lagrange multiplier,
  \item  $x^\star$, $\lambda^\star$, and $\mu^\star$ are locally Lipschitz at~$\bar{\xi}$.
\end{romcases}
\end{theorem}
In particular, note that, $x^\star$ and $\mu^\star$ are feasible on $\calN$, i.e.,  $x^\star(\xi) \in \calX(\xi)$ and $\mu^\star(\xi) \in \R^M_{\geq 0}$ hold for all $\xi \in \calN$.

\section{Quantifying the Lipschitz Continuity of Solution Maps} \label{sec:lipschitz.conitnulity.of.solution.maps}

\cref{thm:sens_continuous} and similar results in \cite{fiaccoSensitivityAnalysisNonlinear1976,robinsonPerturbedKuhnTuckerPoints1974,jittorntrumSolutionPointDifferentiability1984} guarantee the Lipschitz continuity of the solution map $x^\star$, however, they do not quantify the Lipschitz constants of $x^\star$ and $\mu^\star$.

Hence, in this section, we refine these results and give quantitative bounds on the rate of change of $x^\star$. The key insight for this purpose, also used in \cite{fiaccoSensitivityAnalysisNonlinear1976,jittorntrumSolutionPointDifferentiability1984}, is that a KKT point $(x^\star, \mu^\star) \in \calX(\xi) \times \R^m_{\geq 0}$ solves the system
\begin{align} \label{eq:KKT.conditions}
\begin{split}
    F(x^\star, \mu^\star, \xi) := \begin{bmatrix} \grad_x L(x^\star,\mu^\star, \xi) \\ h(x^\star, \xi) \\
      \diag(\mu^\star) g(x^\star, \xi) \end{bmatrix} &= \bbzero \in \R^{n + p + m}.
\end{split}
\end{align}
Under additional assumptions, the implicit function theorem can be applied to \eqref{eq:KKT.conditions}, to express $(x^\star, \mu^\star)$ as a function of $\xi$.

Throughout the remainder of this section, we consider the same setup as in \cref{thm:sens_continuous}. Hence, existence and local Lipschitz continuity of $x^\star: \calN \rightarrow \R^n$ and $\mu^\star: \calN \rightarrow \R^m$ for some neighborhood $\calN \subset \Xi$ around $\bar{\xi}$, for which $(\bar{x}^\star, \bar{\mu}^\star)$ is a regular optimizer, are guaranteed.
Further, we use the shorthand notations $\bfI^\star(\xi) := \bfI^\xi_{x^\star(\xi)}$,  $\bar{\bfI}^\star(\xi) := \bar{\bfI}^\xi_{x^\star(\xi)}$, and, for all $\xi \in \calN$, we define
\begin{align*}
    A(\xi) &\coloneqq \nabla^2_{xx} L(x^\star(\xi), \lambda^\star(\xi), \mu^\star(\xi), \xi))  \in \R^{n\times n}\, ,\\
    B(\xi) &\coloneqq \begin{bmatrix}
        \nabla_x h(x^\star(\xi), \xi) \\
        \nabla_x g_{\bfI^\star(\xi)} (x^\star(\xi), \xi)
        \end{bmatrix}
        \in \R^{(p + | \bfI^\star|) \times n}\, ,\\
    L^\star(\xi)  &:= \nabla^2_{\xi x}  L (x^\star(\xi), \lambda^\star(\xi), \mu^\star(\xi), \xi) \in \R^{n \times r}\, , \\
    G^\star(\xi) &:= \begin{bmatrix}
        \nabla_\xi h(x^\star(\xi), \xi) \\
        \nabla_\xi  g_{\bfI^\star(\xi)} (x^\star(\xi), {\xi})
        \end{bmatrix}
        \in \R^{(p + | \bfI^\star|) \times r}\, .
\end{align*}

If \eqref{eq:basic_opt_prob} does not have any equality constraints and $\bfI^\star(\xi) = \emptyset$, then we follow the convention that $B(\zeta) = 0$ and $G^\star(\xi) = 0$.

In the following, unless there is any ambiguity, we drop the argument from any map whose sole argument is $\xi$.

Finally, we make one simplifying assumption that can, presumably, be slightly relaxed in future work:

\begin{assumption}\label{ass:inv_A}
    In the setup of \cref{thm:sens_continuous}, the matrix $A(\xi)$ is positive definite for all $\xi \in \Xi$.
\end{assumption}

In particular, \cref{ass:inv_A} is satisfied for convex optimization problems with strongly convex objective -- though one can perceive also  scenarios with weaker regularity assumptions. Furthermore, the assumption is common in multi-parametric programming, when the the KKT system needs to be solved explicitly (e.g., \cite{tondelAlgorithmMultiparametricQuadratic2003}). This assumption will be made in forthcoming sections to provide global bounds on the Lipschitz constants of $x^\star$ and~$(\lambda^\star, \mu^\star)$.

To quantify the Lipschitz constant of the solution map we have to distinguish between two cases related to constraints becoming active or inactive when varying $\xi$. Namely, if every active constraint is associated with a positive multiplier (which we refer to as \emph{strict complementarity} below), these constraints remain active for small variations in $\xi$. In this case, the solution maps $x^\star$ and $\mu^\star$ are differentiable. However, when a constraint is active but its multiplier is zero, then differentiability is in general not guaranteed. In this case, we need to make a careful case distinction considering the constraint as being either active or inactive.

\subsection{Sensitivity Analysis under Strict Complementarity}

Given a KKT-point $(x^\star, \lambda^\star, \mu^\star)$ of \eqref{eq:basic_opt_prob} for a given $\xi \in \Xi$, we define the sets of \emph{strongly} active inequality constraints as
\begin{align*}
    \bfIS(x^\star, \mu^\star, \xi) & \coloneqq \{i \in \bfI_{x^\star}^\xi \, | \, \mu^\star_i>0 \}\,
\end{align*}
and we say that $(x^\star, \lambda^\star, \mu^\star)$ satisfies \emph{strict complementary slackness (SCS)} if $\bfIS(x^\star, \mu^\star, \xi)  = \bfI_{x^\star}^\xi$, i.e., if all active inequality constraints are strongly active.

Under strict complementary slackness, the solution maps $x^\star$ and $(\lambda^\star, \mu^\star)$ are continuously differentiable in a neighborhood of $\bar{\xi}$. The following result gives an explicit expression for the respective derivatives.

\begin{proposition}\label{prop:sens_strict_comp}
 Consider the same setup as in \cref{thm:sens_continuous} and let SCS and \cref{ass:inv_A} hold at $(\Bar{x}^\star, \Bar{\lambda}^\star, \Bar{\mu}^\star)$. Then, in a neighborhood $\calN \subset \Xi$ of $\bar{\xi}$, the maps $x^\star: \calN \rightarrow \R^n$ and $(\lambda^\star, \mu^\star): \calN \rightarrow \R^{p} \times \R^m_{\geq0}$ are continuously differentiable with
\begin{equation}\label{eq:opt_deriv}
\begin{split}
    \nabla_\xi x^\star &=
    - \Pi  A^{-1} L^{\star}
    - \Sigma B^\dagger G^{\star}  \in \R^{n \times r}\, , \\ 
     \nabla_\xi\begin{smallbmatrix}
         \lambda^\star \\
        \mu^\star_{\bfI^\star}
    \end{smallbmatrix} &=
     {B^\dagger}^\T A \Sigma \left( A^{-1} {L^\star} - B^\dagger {G^\star} \right) \in \R^{(p + | \bfI^\star|) \times r}\,
\end{split}
\end{equation}
(and $\nabla_\xi \mu_{\bar{\bfI}^\star}^\star = 0$), where $B^\dagger := B^\T ( B B^\T)^{-1} \in \R^{n \times (p + | \bfI^\star |)}$ is the pseudoinverse of $B$, and $\Sigma, \Pi \in \R^{n \times n}$ are given by
\begin{align*}
    \Sigma := A^{-1} B^\T (B A^{-1} B^\T)^{-1} B \quad \text{and} \quad \Pi :=  \bbI - \Sigma \, .
\end{align*}
If \eqref{eq:basic_opt_prob} has only inequality constraints and none of them are active at $\bar{\xi}$, then we have $B^\dagger G^\star = 0$.
\end{proposition}

\begin{proof}
We follow the same procedure as in \cite{fiaccoSensitivityAnalysisNonlinear1976, jittorntrumSequentialAlgorithmsNonlinear1978}. However, \cref{ass:inv_A}, the matrix inversion \cref{lem:block_inverse} in the appendix, and some attention to details allow us to derive the explicit expressions \eqref{eq:opt_deriv} for the derivative of the solution map.

For the moment let $z^\star = (x^\star, \lambda^\star, \mu^\star)$. The implicit function theorem \cite[Ch. 9]{rudinPrinciplesMathematicalAnalysis1976} applied to~\eqref{eq:KKT.conditions} states that $\nabla_{\xi} z^\star = -\left(\nabla_{z} F\right)^{-1} \nabla_\xi F$.
Hence, without loss of generality, assume that the first $| \bfI^\star|$ constraints are (strongly) active and, by SCS, the remaining constraints are inactive. Then, $\nabla F$ can be written as
\begin{align*}
    \nabla_{z} F  & \coloneqq
    \begin{bmatrix}
       A    &                    B^\T       &   \bar{B}^\T \\
      D_\mu B &  \bbzero &\bbzero \\
      \bbzero & \bbzero &  D_{\bar{g}}
    \end{bmatrix} \, , \, \,
    \nabla_{\xi} F& \coloneqq \begin{bmatrix}
    L^\star \\
    D_\mu G^\star \\
    \bbzero
    \end{bmatrix} \, ,
\end{align*}
where we partitioned and regrouped the matrices according to $(x^\star, [\lambda^\star,  \mu_{\bfI^\star}^\star] ,\mu_{\bar\bfI^\star}^\star)$, and used
$\bar{B}(\xi) \coloneqq \nabla_x g_{\bar{\bfI}^\star(\xi)} (x^\star(\xi), \xi)$. Further, we have used
$D_\mu :=\diag( [\bbone_p^\T \,  \mu_{\bfI^\star}^\T ])$ and $D_{\bar{g}} := \diag(g_{\bar{\bfI}^\star})$.
Denote the upper left part of $\nabla_{z} F$ by
\begin{align}
     M \coloneqq   \begin{bmatrix}
    A && B^\T \\ D_{\mu} B&& \bbzero
    \end{bmatrix} \, .
\end{align}
Then, using \cref{lem:block_inverse} (and assuming for the moment that $M$ is invertible), we have that
\begin{align*}
    \left( \nabla_{(x, \mu)} F \right)^{-1} = \begin{bmatrix}
        M^{-1} & \star \\ \bbzero &  \star
    \end{bmatrix}   \, ,
\end{align*}
where $\star$ denotes non-zero, but irrelevant components.

Hence, we can already conclude from the expression for $\nabla_\xi z^\star$ that $\nabla_\xi \mu_{\bar{\bfI}^\star}^\star = 0$ for the inactive constraints.

Note that $A$ is invertible by \cref{ass:inv_A}, and the Schur complement of $M$ given by $- D_\mu B A^{-1} B^T$ is invertible because $D_\mu$ in invertible by SCS. Hence, $M$ is nonsingular, and its inverse can be be computed according to \cref{lem:block_inverse}
as
\begin{align*}
     M^{-1} = \begin{bmatrix} M_1 && M_2 \\ M_3 && M_4 \end{bmatrix},
\end{align*}
where
\begin{align*}
    M_1
     &= \left( \bbI - A^{-1} B^\T( D_\mu B A^{-1} B^\T)^{-1} D_\mu B \right) A^{-1} \\
     &= (\bbI - A^{-1} B^\T( B A^{-1} B^\T)^{-1} B) A^{-1} = \Pi A^{-1} \, , \\
    M_2
     &= A^{-1} B^\T (D_\mu B A^{-1} B^\T)^{-1} = \\
     &= A^{-1} B^\T (B A^{-1} B^\T)^{-1} D_\mu^{-1} \\
     &= A^{-1} B^\T (B A^{-1} B^\T)^{-1} D_\mu^{-1} \underbrace{D_\mu B B^\dagger D^{-1}_\mu}_{= \bbI} =\Sigma B^\dagger D_\mu^{-1} \, , \\
    M_3
     &= ( D_\mu  B A^{-1} B^\T)^{-1} D_\mu B A^{-1} \\
     &= \underbrace{{B^\dagger}^\T A A^{-1} B^\T }_{\bbI} ( BA^{-1}B^\T)^{-1} B A^{-1} = {B^\dagger}^\T A \Sigma A^{-1} \, ,\\
    M_4
     &= -(D_\mu BA^{-1}B^\T)^{-1} = -(BA^{-1}B^\T)^{-1} D_\mu^{-1} \\
     &= - {B^\dagger}^\T A A^{-1} B^\T  (BA^{-1}B^\T)^{-1} D_\mu^{-1} D_\mu B B^\dagger D^{-1}_\mu\\
     &= - { B^\dagger}^\T A \Sigma B^\dagger D_\mu^{-1}.
\end{align*}

Finally, recalling $\nabla_{\xi} z^\star = -\left(\nabla_{z} F\right)^{-1} \nabla_\xi F$, we get
\begin{align*}
    \begin{bmatrix} \nabla_\xi x^\star \\
    \nabla_\xi \lambda^\star \\
    \nabla_\xi \mu^\star_{\bfI^\star} \end{bmatrix} = -
\begin{bmatrix}    M^{-1} && \star \end{bmatrix} \nabla_\xi F
\end{align*}
which yields the desired expressions.
\end{proof}

We defer an in-depth discussion of \eqref{eq:opt_deriv} to \cref{sec:disc}, where we explore special cases.
One special feature of \eqref{eq:opt_deriv} that we require though is the fact that $\Pi$ and $\Sigma$ can be interpreted as (oblique) projections onto the kernel of $B$ and its orthogonal complement, respectively. Consequently, exploiting non-expansivity, their norm can be bounded independently of~$B$:

\begin{lemma}\label{lem:proj_bound}
    Consider the setup of \cref{prop:sens_strict_comp}. For all $\xi \in \calN$, it holds that $\max\{\| \Pi(\xi) \| , \|\Sigma(\xi) \| \} \leq  \sqrt{\lambda^{\max}_{A(\xi)}/\lambda^{\min}_{A(\xi)}}$.
\end{lemma}

\begin{proof}
    As before, we drop the argument $\xi$. Consider
    \begin{align}\label{eq:proj_prob}
        \underset{w}{\minimize} \, \tfrac{1}{2} (v - w)^\T A (v - w) \, \, \st \, B w = 0 \,.
    \end{align}
    Since $A$ is positive definite for all $\xi \in \Xi$ by \cref{ass:inv_A},~\eqref{eq:proj_prob} can be rewritten as
    \begin{align}\label{eq:proj_prob_reform}
        \underset{\tilde{w}}{\minimize} \quad \tfrac{1}{2} \| Q v - \tilde{w} \|^2 \quad \st \quad \tilde{B} \tilde{w} = 0
    \end{align}
    where $Q := A^{1/2}$ is the unique symmetric positive definite matrix such that $A = Q^2$, $\tilde{w} := Q w$, and $\tilde{B} := B Q^{-1}$.

    Since~\eqref{eq:proj_prob_reform} is a Euclidean projection onto the kernel of $\tilde{B}$, its solution can be stated explicitly as \cite[eq. 5.13.3]{meyerMatrixanalysisapplied2000}
    \begin{align*}
        \tilde{w}^\star
        & = Q v - \tilde{B}^T ( \tilde{B} \tilde{B}^\T )^{-1} \tilde{B} Q v \\
        & = Q( \bbI -  A^{-1} B^\T ( B A^{-1} B^\T )^{-1} B ) v  = Q \Pi v \, .
    \end{align*}
    By non-expansivity of an orthogonal projection onto a subspace, we have $\| \tilde{w}^\star \| \leq \| Q v\|$, and therefore we can write
    \begin{align*}
        \lambda^{\min}_Q \| \Pi v \| \leq \| Q \Pi v \| = \| \tilde{w}^\star \| \leq \| Q v \| \leq \lambda^{\max}_Q \| v \| \, .
    \end{align*}
    Rearranging the leftmost and rightmost terms together with the fact that $\lambda^{\max}_A = (\lambda^{\max}_Q)^2$ and  $\lambda^{\min}_A = (\lambda^{\min}_Q)^2$ yields the desired bound for $\|\Pi \|$.

    The projection on the orthogonal complement of $\ker \tilde{B}$ is given by $\Sigma := \bbI - \Pi$ \cite[eq. 5.13.6]{meyerMatrixanalysisapplied2000}. The same bound as for $\|\Pi\|$ holds for~$\| \Sigma \|$ since $\Sigma$ is also non-expansive.
\end{proof}

\cref{lem:proj_bound} allows us to bound~\eqref{eq:opt_deriv} and thereby establish bounds for the local Lipschitz constants of $x^\star$ and $\mu^\star$ at $\bar{\xi}$.

\begin{corollary}\label{cor:lip_strict_comp}
   Consider the setup of \cref{prop:sens_strict_comp}. Then,
   \begin{align*}
      \left\|\nabla_\xi x^\star(\bar{\xi}) \right \| \leq \ell_{x^\star} \quad \text{and} \quad
       \left\|\nabla_\xi \begin{smallbmatrix}
         \lambda^\star \\ \mu^\star
    \end{smallbmatrix} \right\| \leq \ell_{(\lambda^\star,\mu^\star)}
   \end{align*}
   holds for all $\xi \in \calN$ where
   \begin{align}\label{eq:lip_bound_strict_comp}
   \begin{split}
     \ell_{x^\star}
        & :=  \sqrt{\frac{\lambda^{\max}_{A(\bar{\xi})}}{\lambda^{\min}_{A(\bar{\xi})}}} \left(
             \frac{\| L^\star(\bar{\xi}) \|}{ \lambda^{\min}_{A(\bar{\xi})}}
          + \frac{\| G^\star(\bar{\xi})\|}{\sigma^{\min}_{B^\T(\bar{\xi})}}
         \right) \, , \\
       \ell_{(\lambda^\star,\mu^\star)}
        &:= \frac{{\lambda^{\max}_{A(\bar{\xi})}}^{3/2}}{\sigma^{\min}_{B^\T(\bar{\xi})} {\lambda^{\min}_{A(\bar{\xi})}}^{1/2}} \left(
             \frac{\| L^\star(\bar{\xi}) \|}{ \lambda^{\min}_{A(\bar{\xi})}}
          + \frac{\| G^\star(\bar{\xi})\|}{\sigma^{\min}_{B^\T(\bar{\xi})}}
         \right) \,.
    \end{split}
   \end{align}
\end{corollary}
\begin{proof}
 Using singular value decomposition $B^\T=V^\T \Lambda^\T U$, where $B^\T$ (i.e., $\Lambda^\T$) has a full rank since the columns of $B$ are linearly independent. Next we note that $B^\dagger=V^\T \Lambda^\dagger U$ also has a full rank and corresponds to the right pseudoinverse of $B$, where $\Lambda^\dagger$ is the right pseudoinverse of $\Lambda$. Comparing the expressions for $B^\T$ and $B^\dagger$ we have that for the minimal and maximal singular values $\sigma^{\max}_{B^\dagger}=1 / \sigma^{\min}_{B^\T}$ and $\sigma^{\min}_{B^\dagger}=1 / \sigma^{\max}_{B^\T}$ hold. Consequently,  $\norm{B^\dagger} = 1 / \sigma^{\min}_{B^\T(\bar{\xi})}$ holds. Finally, the bounds in \eqref{eq:lip_bound_strict_comp} follows from applying the triangle inequality, Cauchy-Schwarz, and \cref{lem:proj_bound} to~\eqref{eq:opt_deriv}.
\end{proof}
\subsection{Sensitivity Analysis for General Regular Optimizers}
If the SCS assumption in \cref{prop:sens_strict_comp} does not hold, then at least one constraint $i$ is \emph{weakly active}, i.e., $g_i(x^\star(\xi))=0$ but $\mu_i^\star(\xi)=0$. Weakly active constraints occur, for example, when the unconstrained minimizer of $f$ lies exactly on a constraint surface. In this case, the constraint is active, but is not ``pushing'' the optimizer inwards.

Whenever we vary $\xi$ and a constraint changes from being strongly active to inactive (or vice versa), the constraint is momentarily weakly active for some $\xi$. This is a consequence of the continuity of $(x^\star, \mu^\star)$.

For values of $\xi$, for which some constraints are weakly active, $(x^\star, \lambda^\star, \mu^\star)$ might not be differentiable, i.e., \eqref{eq:opt_deriv} might not apply.\footnote{However, it is generally possible to compute directional derivatives of $(x^\star, \lambda^\star, \mu^\star)$ at any point in any direction $\Delta \xi$.} Consequently, the best we can hope for in the general case where SCS does not hold, is a bound on the Lipschitz constant of $(x^\star, \lambda^\star, \mu^\star)$.

To deal with weakly active constraints, we follow the approach in~\cite{jittorntrumSequentialAlgorithmsNonlinear1978}. For this purpose, consider the same setup as in \cref{prop:sens_strict_comp} and let $(\bar{x}^\star, \bar{\lambda}^\star, \bar{\mu}^\star)$ be a regular optimizer for~\cref{eq:basic_opt_prob} for $\bar{\xi} \in \Xi$, satisfying LICQ at $\bar{x}^\star$, KKT, and SSOSC but not necessarily SCS.
Further, let $\bfR$ be any index set satisfying
\begin{equation}\label{eq:idx_set_R}
\begin{aligned}
    \bfIS(\bar{x}^\star, \bar{\mu}^\star, \bar{\xi}) \subset \bfR \subset \bfI^\star(\bar{\xi}) \, ,
\end{aligned}
\end{equation}
that is, $\bfR$ includes all strongly active constraints and an arbitrary subset of weakly active constraints at $(\bar{x}^\star, \bar{\lambda}^\star, \bar{\mu}^\star, \bar{\xi})$. Next, consider the equality-constrained problem
\begin{equation}\label{eq:eq_cstr_prob}
\begin{split}
    \minimize \quad & f(x, \xi) \\
          \st \quad & h(x, \xi) = 0 \\
                    & g_\bfR(x, \xi) = 0 \, .
\end{split}
\end{equation}
By considering \eqref{eq:eq_cstr_prob} we make a choice for every constraint that is weakly active at $(\bar{x}^\star, \bar{\lambda}^\star, \bar{\mu}^\star)$ to consider it like a strongly active constraint (by including it in $\bfR$) or to ignore it.

Note that \cref{prop:sens_strict_comp,cor:lip_strict_comp} can be applied to~\eqref{eq:eq_cstr_prob}. In particular, the SCS condition is vacuous because \eqref{eq:eq_cstr_prob} does not have any inequality constraints.
Hence, let $x^\star_\bfR : \calN_{\bfR} \rightarrow \R^n$ and $(\lambda^\star_\bfR, \mu^\star_\bfR) : \calN_{\bfR} \rightarrow \R^{(p + \abs{\bfR})}$ denote the primal and dual solution map of \eqref{eq:eq_cstr_prob} in a neighborhood $\calN_\bfR$ of $\bar{\xi}$.

In the proof of \cite[Thm 2.3.2]{jittorntrumSequentialAlgorithmsNonlinear1978} it was established\footnote{To establish this equivalence the author passes through an exact penalty reformulation, which is beyond the scope of this paper.} that for every $\xi$ in a neighborhood $\calN \subset \bigcap_\bfR \calN_\bfR$ of $\bar{\xi}$ the following holds: If $(x^\star(\xi), \lambda^\star(\xi), \mu^\star(\xi))$ is a regular optimizer of \eqref{eq:basic_opt_prob}, then there exists an index set $\bfR$ satisfying~\eqref{eq:idx_set_R} such that $(x^\star(\xi), \lambda^\star(\xi), \mu^\star_\bfR(\xi))$ is a regular optimizer of \eqref{eq:eq_cstr_prob}.
In particular, we have $x^\star(\xi) = x^\star_\bfR(\xi)$, $\lambda^\star(\xi) = \lambda^\star_\bfR(\xi)$, and $\mu^\star_i(\xi) = 0$ for all $i \notin \bfR$.

Note that the set $\bfR$ for which this equivalence of solutions between \eqref{eq:basic_opt_prob} and \eqref{eq:eq_cstr_prob} holds, depends on $\xi$ and there does in general not exist a single set $\bfR$ that works for the entire neighborhood $\calN$ of $\bar{\xi}$.

This key insight can be used to establish bounds on the Lipschitz constants of $x^\star$ and $\mu^\star$ at $\bar{\xi}$ in the absence of SCS:

\begin{theorem}\label{thm:sens_degen}
  Consider the same setup as in \cref{thm:sens_continuous} and let \cref{ass:inv_A} hold at $(\Bar{x}^\star, \Bar{\lambda}^\star, \Bar{\mu}^\star)$. Then, $x^\star: \calN \rightarrow \R^n$ and $(\lambda^\star, \mu^\star): \calN \rightarrow \R^{p + m}$ are locally Lipschitz at $\bar{\xi}$ with bounds for the Lipschitz constants given by $\ell_{x^\star}$ and  $\ell_{(\lambda^\star,\mu^\star)}$ as in~\eqref{eq:lip_bound_strict_comp}.
\end{theorem}

The difference in \eqref{eq:lip_bound_strict_comp} for \cref{cor:lip_strict_comp} and \cref{thm:sens_degen} lies in the fact that for \cref{cor:lip_strict_comp} the set $\bfI^\star$ (required for $A, B, L^\star$ and $G^\star$) contains only strictly active constraints, whereas for \cref{thm:sens_degen} it may also contain weakly active constraints.

\begin{proof}
     Given $\bfR$ satisfying \eqref{eq:idx_set_R}, let $x^\star_\bfR: \calN \rightarrow \R^n$, $\lambda^\star_\bfR: \calN \rightarrow \R^m$, and $\mu^\star_\bfR: \calN \rightarrow \R^{\abs{\bfR}}$ the solution maps of \eqref{eq:eq_cstr_prob} in a neighborhood $\calN$ of $\bar{\xi}$.
    Consequently, we apply \cref{cor:lip_strict_comp} to \eqref{eq:eq_cstr_prob} for every $\bfR$ satisfying \eqref{eq:idx_set_R}, and we let $\ell^\bfR_{x^\star}$ and $\ell^\bfR_{{(\lambda^\star,\mu^\star)}}$ denote the respective bounds on the Lipschitz constant of $x^\star_\bfR$ and $(\lambda_{\bfR}^\star, \mu^\star_{\bfR})$ at $\bar{\xi}$ according to \eqref{eq:lip_bound_strict_comp}. Similarly, let $A_\bfR(\bar{\xi})$, $B_\bfR(\bar{\xi})$, $L^\star_\bfR(\bar{\xi})$ and $G^\star_\bfR(\bar{\xi})$ denote the corresponding quantities for~\eqref{eq:eq_cstr_prob}, evaluated at $\bar{\xi}$.

    Since, for every $\xi \in \Xi$ the solution $(x^\star(\xi),\mu^\star(\xi))$ is equal to $(x_\bfR^\star(\xi),\mu^\star(\xi))$ for some $\bfR$ satisfying \eqref{eq:idx_set_R}, we can upper bound the local Lipschitz constants of $x^\star$ and $(\lambda^\star, \mu^\star)$ at $\bar{\xi}$ by maximizing over all (finite) possibilities of $\bfR$ satisfying \eqref{eq:idx_set_R}. Namely, for $x^\star$ we have
    \begin{align*}
        \sup_{\xi \rightarrow \bar{\xi}} \tfrac{\| x^\star(\xi) - x^\star(\bar{\xi}) \|}{\| \xi - \bar{\xi} \|}
        \leq \max_{\bfR} \sup_{\xi \rightarrow \bar{\xi}} \tfrac{\| x^\star_\bfR(\xi) - x^\star(\bar{\xi}) \|}{\| \xi - \bar{\xi} \|} \leq \max_{\bfR} \ell_{x^\star}^\bfR \, ,
    \end{align*}
    and the case for $(\lambda^\star, \mu^\star)$ follows analogously.

    In particular, we claim that this maximum is achieved exactly for the choice $\bfR = \bfI^\star(\bar{\xi})$ where all strongly and weakly active constraints are considered.

    To see this, note the following: Since lifting a vector to higher dimensions by adding components can only increase its 2-norm, we have $\| G^\star_\bfR (\bar{\xi}) \| \leq \| G^\star_{\bfI^\star(\bar{\xi})} (\bar{\xi}) \|$ for any $\bfR$ satisfying \eqref{eq:idx_set_R}. Analogously, we have $\sigma^{\min}_{B^\T_\bfR(\bar{\xi})} \geq \sigma^{\min}_{B^\T_{\bfI^\star (\bar{\xi})}}$ since adding columns to $B^\T$ can only reduce its minimum singular value.

    Next, at $\bar{\xi}$ we that $\mu^\star(\bar{\xi}) = \mu^\star_\bfR(\bar{\xi})$ for any $\bfR$ satisfying \eqref{eq:idx_set_R}. Further, since $\mu^\star_{i}(\bar{\xi}) = 0$, for all weakly active constraints with index $i$, we have $ L^\star (\xi) = L^\star_\bfR (\xi) = L^\star_{\bfI^\star(\bar{\xi})} (\xi) $ and $ A(\bar{\xi}) = A_\bfR(\bar{\xi}) =  A_{\bfI^\star\textbf{}(\bar{\xi})}(\bar{\xi})$. In other words, at $\bar{\xi}$ the weakly active constraints do not affect the quantities $A$ and $L^\star$ which are derivatives of the Lagrangian.

    Applying these facts to \eqref{eq:lip_bound_strict_comp}, it follows that $\ell^\bfR_{x^\star} \leq \ell_{x^\star}^{\bfI^\star(\bar{\xi})}$ and  $\ell^\bfR_{(\lambda^\star,\mu^\star)} \leq \ell_{{(\lambda^\star,\mu^\star)}}^{\bfI^\star(\bar{\xi})}$, for all $\bfR$ satisfying \eqref{eq:idx_set_R} and thus we have
       $\ell_{x^\star} = \max_{\bfR} \ell^{\bfR}_{x^\star} =  \ell^{\bfI^\star(\bar{\xi})}_{x^\star}$ and
         $\ell_{{(\lambda^\star,\mu^\star)}} =  \max_{\bfR} \ell^{\bfR}_{{(\lambda^\star,\mu^\star)}} \!=\!  \ell^{\bfI^\star(\bar{\xi})}_{{(\lambda^\star,\mu^\star)}} $ which completes the proof.
\end{proof}
\section{Assumptions for Global Solution Maps} \label{sec:assumptions.for.global.solution.maps}
\cref{thm:sens_degen} provides bounds on the local Lipschitz constants for $x^\star$ and $\mu^\star$ at $\bar{\xi}$. Next, we provide a set of assumptions to give global bounds on the sensitivity. Naturally, some of these assumptions are restrictive, but they provide important intuition and use cases for the general result in~\cref{thm:sens_degen}.

For simplicity, instead of \eqref{eq:basic_opt_prob}, we focus on problems with only inequality constraints which require a more careful investigation than equality constraints as shown in the previous analysis, i.e., we consider
\begin{align}\label{eq:basic_opt_prob_ineq}
    \begin{split}
        \underset{x}{\minimize} \quad & f(x, \xi)\\
        \st \quad & g(x, \xi) \leq 0 \, ,
      \end{split}
\end{align}
for the remainder of this section. All of the following statements and assumptions can be generalized.

First, it is necessary that~\eqref{eq:basic_opt_prob_ineq} admits a solution for every~$\xi$:
\begin{assumption}\label{ass:feasible}
    The problem \eqref{eq:basic_opt_prob_ineq} is feasible for all $\xi \in \Xi$.
\end{assumption}

Next, in order to guarantee that the solution map $x^\star: \Xi \rightarrow \R^n$ is single-valued for all $\xi \in \Xi$, it is convenient, if not necessary, to assume (strong) convexity:
\begin{assumption}\label{ass:cvx}
    For~\eqref{eq:basic_opt_prob_ineq}, let $f$ be  $\alpha$-strongly convex and have a $\beta$-Lipschitz gradient $\nabla_x f$ in $x$ for all $\xi \in \Xi$. Further, for all $i = 1, \ldots, m$, let $g_i$ be convex and $\nabla_x g$ be $\ell_i$-Lipschitz in $x$ for all $\xi \in \Xi$.
\end{assumption}
\cref{ass:cvx} not only guarantees that \eqref{eq:basic_opt_prob} admits a unique optimizer for every $\xi \in \Xi$, it also implies that any KKT point of \eqref{eq:basic_opt_prob_ineq} (trivially) satisfies the SSOSC since $\nabla^2_{xx} L$ is positive definite for all $(x, \mu, \xi) \in \R^n \times \R^m_{\geq 0} \times \Xi$. \cref{ass:cvx} also implies the lower bound $\lambda^{\min}(A(\xi)) \geq \alpha$ for all $\xi \in \Xi$ and thus replaces \cref{ass:inv_A}.
 Lipschitz continuity of $\nabla_x f$ and $\nabla_x g$ is required to upper bound $\lambda^{\max}(A)$ as discussed below.
Moreover, we generally require the active constraints to have uniformly full rank:
\begin{assumption}\label{ass:unif_licq}
    Given~\eqref{eq:basic_opt_prob_ineq}, there exists $\omega > 0$ such that, for all $\xi \in \Xi$ and all $x \in \calX(\xi)$ such that $\bfI_x^\xi \neq \emptyset$, we have
    \begin{align*}
        \omega^2 \bbI \preceq \nabla_{x} g_{\bfI^\xi_{x}}(x, \xi) \nabla_{x} g_{\bfI^\xi_{x}}(x, \xi)^\T \, .
    \end{align*}
\end{assumption}
In particular, \cref{ass:unif_licq} guarantees that LICQ is satisfied at all $x \in \calX(\xi)$, for all $\xi \in \Xi$, and that $\sigma^{\min}_{B^\T(\xi)} \geq \omega$ in~\eqref{eq:lip_bound_strict_comp} is lower bounded away from zero whenever at least one constraint is active. Note that \cref{ass:unif_licq} is independent of the choice of the cost function.

Further note that, under \cref{ass:cvx,ass:unif_licq}, \eqref{eq:basic_opt_prob_ineq} has a unique regular optimizer for every $\xi \in \Xi$. Hence, by \cref{thm:sens_continuous}, continuous solution (and multiplier) maps exist around every $\xi \in \Xi$. This implies that the maps $x^\star: \Xi \rightarrow \R^n$ and $\mu^\star: \Xi \rightarrow \R^m$ exist globally on all of~$\Xi$ and are continuous.

Finally, because $A = \nabla^2_{xx} f + \sum_i \mu_i \nabla_{xx} g$ and $L^\star = \nabla_{\xi x} f + \sum_i \mu_i \nabla^2_{\xi x} g$ are weighted sums over $\mu$, we generally require an upper bound on the dual multiplier $\mu^\star$.

\begin{assumption}\label{ass:dual_bound}
    There exist $\zeta_i > 0$ such that, for all $\xi \in \Xi$ and every KKT point $(x^\star(\xi), \mu^\star(\xi))$ of~\eqref{eq:basic_opt_prob_ineq}, we have $ \mu^\star_i(\xi) \leq \zeta_i$ for all $i=1, \ldots, m$.
\end{assumption}

Finding (tight) upper bounds on the dual multipliers is tricky and depends on the problem structure.

If $\| \nabla_x f(x, \xi) \| \leq B_f $ and $\|\nabla_x g_i(x, \xi)\| \leq B_{g_i} $ are uniformly bounded for all $\xi \in \Xi$, all $x \in \calX(\xi)$ and all $i = 1, \dots, m$, then it is easy to see that $\mu_i^\star \leq \zeta_i := B_f / B_{g_i}$.

Another possibility, documented in \cite[p.647]{bertsekasdescentnumericalmethod1973}, applies specifically to convex optimization problems (i.e., under \cref{ass:cvx}). Namely, if for all $\xi \in \Xi$ a strict lower bound $\tilde{f}^\star_\xi$ (i.e., $\tilde{f}^\star_\xi < f(x, \xi)$ for all $x \in \calX(\xi)$) and a strictly feasible point $\tilde{x}_\xi$ (i.e., $g(\tilde{x}_{\xi}, \xi) < 0$) are known, then we can choose
\begin{align*}
    \zeta_i = \sup_{\xi \in \Xi} \frac{f(\tilde{x}_\xi, \xi) - \tilde{f}^\star_\xi}{- g_i(\tilde{x}_\xi)} \, .
\end{align*}
Combining \cref{ass:cvx,ass:dual_bound} we can guarantee that
\begin{align*}
    A  = \nabla_{xx}^2 f    + \sum\nolimits_{i=1}^m \mu_i^\star \nabla^2_{xx} g \preceq \left(\beta + \sum\nolimits_i \zeta_i \ell_i \right) \bbI\, .
\end{align*}
In other words, we have  $\lambda^{\max}(A(\xi)) \leq \beta + \sum_i \zeta_i \ell_i $.

The assumptions made so far can be summarized in the following statement:

\begin{corollary} \label{cor:ell.lipschitz.optimizer}
    Under \cref{ass:feasible,ass:cvx,ass:unif_licq,ass:dual_bound,} the primal and dual solution maps $x^\star: \Xi \rightarrow \R^n$ and $\mu^\star: \Xi \rightarrow \R^m $ of \eqref{eq:basic_opt_prob_ineq} are single-valued and Lipschitz continuous (on $\Xi$) with respective Lipschitz constants
\begin{align}
\begin{split} \label{eq:ell.lipschitz.optimizer}
    \ell_{x^\star} &:= \sqrt{\frac{\beta + \sum_i \zeta_i \ell_i}{\alpha}} \left( \frac{\bar{L}^\star}{\alpha} + \frac{\bar{G}^\star}{\omega} \right) \, , \\
    \ell_{\mu^\star} &:= \frac{(\beta + \sum_i \zeta_i \ell_i)^{3/2}}{\alpha^{1/2} \omega} \left( \frac{\bar{L}^\star}{\alpha} + \frac{\bar{G}^\star}{\omega} \right) \, ,
    \end{split}
\end{align}
where $\bar{L}^\star :=\underset{\xi \in \Xi} {\sup} \| L^\star(\xi) \|$ and $\bar{G}^\star := \underset{\xi \in \Xi} {\sup} \| G^\star(\xi) \|$ are upper bounds on $L^\star$ and $G^\star$ over $\Xi$, respectively.
\end{corollary}

It remains to establish bounds on $L^\star$ and $G^\star$. Such bounds, however, are highly problem dependent, and therefore we consider only special cases in the next section.

\section{General Discussion of Sensitivity Bounds}\label{sec:disc}

We now discuss special cases for the established sensitivity bounds and how to relax the differentiability of $f$ and $g$.

\begin{figure}[bt]
	\centering
	\subfloat[\label{fig:opt_traj2}]{
		\centering
		\includegraphics[width=.35\columnwidth]{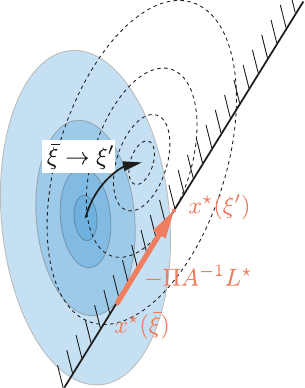}
	}
	\hspace{.5cm}
	\subfloat[\label{fig:opt_traj1}]{
		\centering
		\includegraphics[width=.35\columnwidth]{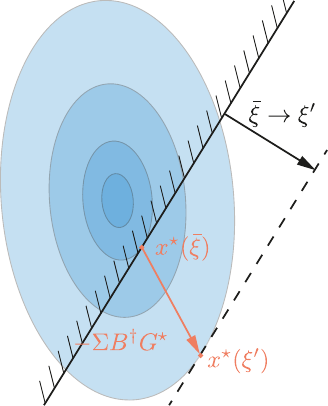}
		}
	\caption{Evolution of the optimizer for (a) stationary constraints and perturbed cost function; (b) fixed objective and linear constraint with righthand side perturbation.}
\end{figure}

\subsection{Special Cases}
\subsubsection*{Stationary Constraints}
We first assume that the constraints are stationary and only the cost function is perturbed, i.e., we consider the problem
\begin{align}
    \minimize \, \, f(x, \xi) \quad \st \, \, g(x) \leq 0
\end{align}
and assume that \cref{ass:feasible,ass:cvx,ass:unif_licq,ass:dual_bound} are satisfied. Let $(\bar{x}^\star, \bar{\mu}^\star)$ be a regular optimizer for $\bar{\xi}$ for which SCS holds, and therefore \cref{prop:sens_strict_comp} is applicable.

In this case, $G^\star(\xi)$ is vacuous and we have $L^\star = \nabla_{x\xi}^2 L = \nabla^2_{x \xi} f$. Hence, the evolution of $x^\star$ and $\mu^\star$ around $\bar{\xi}$ is governed by the first terms in \eqref{eq:opt_deriv}. In particular, we have $\nabla_\xi x^\star = - \Pi A^{-1} {L^\star}$, which is independent of $B^\dagger$ (and hence of $\sigma^{\min}_{B^\T}$). In other words, the conditioning of the constraints (as quantified by \cref{ass:unif_licq}) does not affect the sensitivity of $x^\star$ (however, it does affect the sensitivity of $\mu^\star$).

In this context, we can also make sense of the operator $\Pi$ which projects the quantity $A^{-1} {L^\star}$ onto $\ker B^\T$ which is the tangent space of the surface of (strongly) active constraints. In other words, $\Pi$ guarantees that $x^\star$ can only change along the surface spanned by the (strongly) active constraints, which is plausible since the constraints are stationary. \cref{fig:opt_traj2} illustrates this effect.

If we further assume that no constraints are active at $\bar{x}^\star$, we find ourselves in a situation of unconstrained optimization. In this case, \eqref{eq:opt_deriv} reduces to $\nabla_\xi x^\star = - \left(\nabla^2_{xx}  f(\xi) \right)^{-1} \nabla^2_{\xi x} f$. This results can also be derived by a direct application of the implicit function theorem, for which $\ell_{x^\star}$ reduces to
\begin{align}\label{eq:unconstrained.ell.constant}
    \ell_{x^\star}\coloneqq \frac{\ell_f}{\alpha}, \quad \text{where} \quad \ell_f \coloneqq \underset{\xi' \in \Xi}{\sup} \norm{\nabla^2_{\xi x} f(x^\star(\xi'),\xi')} \, .
\end{align}
In particular, for unconstrained time-varying optimization, the assumption that $\nabla_x f$ is $\beta$-Lipschitz is not required.

\subsubsection*{Translational Perturbation of Objective}

A particular class of objective functions is when the perturbation is in the form of a translation, i.e., we consider
\begin{align}\label{eq:prob_translation_pert}
    \minimize \, \, &\hat{f}(x - c(\xi)) \quad   \st \, \,  g(x) \leq 0  \, ,
\end{align}
where $\hat{f}: \R^n \rightarrow \R$ is $\alpha$-strongly convex and $\beta$-smooth and $c : \R^r \rightarrow \R^n$ is continuously differentiable and $\ell_c$-Lipschitz. Otherwise, let \cref{ass:feasible,ass:cvx,ass:unif_licq,ass:dual_bound} hold.
Considering $f(x,\xi)=\hat{f} (x-c(\xi))$, the structure of \eqref{eq:prob_translation_pert} implies that $\nabla_x f(x, \xi) = \nabla \hat{f}(x - c(\xi))$, and therefore $\nabla^2_{\xi x} f(x, \xi) = - \nabla^2_{xx} \hat{f}(x - c(\xi)) \nabla c(\xi)$, and therefore $\| L^\star\| \leq \beta \ell_c$ (since $\nabla_{\xi x} g = 0$). In fact, $c$ does not need to be continuously differentiable, but only Lipschitz as discussed in the forthcoming \cref{sec:relax_diff}.

\subsubsection*{Right-hand Constraint Perturbations}

Next, we consider a stationary cost function and a constraint parametrization that takes the form of a right-hand side perturbation, i.e.,
\begin{align}\label{eq:prob_bound_pert}
    \minimize \, \, f(x) \quad \st \, \, u(x) \leq v(\xi) \, ,
\end{align}
where $u: \R^n \rightarrow \R^m$ and $v: \Xi \rightarrow \R^m$ are twice continuously differentiable, and
\cref{ass:feasible,ass:cvx,ass:unif_licq,ass:dual_bound} hold.

Let $(\bar{x}^\star, \bar{\mu}^\star)$ be a regular optimizer for $\bar{\xi}$ and let SCS hold, such that $x^\star$ and $\mu^\star$ are differentiable at $\bar{\xi}$.

Exploiting the structure of \eqref{eq:prob_bound_pert}, we have $\nabla^2_{\xi x} g \equiv 0$, and $\nabla ^2_{\xi x} f \equiv 0$, and thus $L^\star \equiv 0$. Consequently, according to \cref{prop:sens_strict_comp} it holds that $\nabla_\xi x^\star = - \Sigma B^\dagger {G^\star}$.
Recall from \cref{lem:proj_bound} and its proof, that $\Sigma$ denotes an orthogonal projection onto the space spanned by $ Q^{-1} \nabla g_i (\xi)$ for all $i \in \bfI^\star(\xi)$ where $Q := A^{1/2}$. Roughly speaking, as the constraints vary according to~\eqref{eq:prob_bound_pert}, the optimizer can only move $Q^{-1}$-orthogonally to the surface spanned by the active constraints. This behavior is illustrated in \cref{fig:opt_traj1}.

\subsubsection*{Linear Constraints}

Yet, a more special case is to assume that $\calX(\xi)$ is a polyhedron, that is, to consider
\begin{align}\label{eq:prob_bound_pert2}
    \minimize \, \, f(x) \quad \st \, \, U x \leq v(\xi) \, ,
\end{align}
where $U \in \R^{m \times n}$ is assumed to have full row rank (thus satisfying \cref{ass:unif_licq}). The fact that in this case we have that $\nabla_{xx}^2 g = 0$ (and that $L^\star = \nabla^2_{\xi x} f$) obviates the need for \cref{ass:dual_bound}, i.e., there is no need to estimate the upper bound of the dual multipliers. Hence, the bounds $\ell_{x^\star}$ and $\ell_{\mu^\star}$,  in \cref{cor:ell.lipschitz.optimizer} simplify to
\begin{align}\label{eq:constant.prob_bound_pert2}
    \ell_{x^\star} := \sqrt{\frac{\beta}{\alpha}} \frac{\ell_v}{\omega} \quad \text{and} \quad
    \ell_{\mu^\star} := \frac{\beta^{3/2}}{\alpha^{1/2}} \frac{\ell_v}{\omega^2} \,,
\end{align}
where $\ell_v := \underset{\xi \in \Xi}{\sup} \| \nabla v(\xi) \|$ is the Lipschitz constant of~$v$.

\subsection{Relaxation of Differentiability Assumption}\label{sec:relax_diff}
All results thus far require $f$ and $g$, from \eqref{eq:basic_opt_prob_ineq}, to be twice continuously differentiable in $(x,\xi)$. Whether it is possible to directly relax this (possibly restrictive) assumption without losing Lipschitz continuity is still an open question.

However, this issue sometimes can be circumvented by a perturbation mapping that isolates the non-smoothness.
Namely, assume that there exist a map $\tilde{\xi}: \R^p \rightarrow \R^r$ that is $\ell_{\tilde{\xi}}$-Lipschitz. Then, based on \cref{eq:basic_opt_prob_ineq}, we consider the problem
\begin{align}\label{eq:time_varying_prob}
    \begin{split}
        \underset{x}{\minimize} \quad & f(x, \tilde{\xi}(\psi)),\\
        \st \quad & g(x, \tilde{\xi}(\psi)) \leq 0 \,
      \end{split}
\end{align}
parametrized in $\psi$. If \eqref{eq:basic_opt_prob_ineq} admits a $\ell_{x^\star}$-Lipschitz solution map with respect to $\xi$, then \eqref{eq:time_varying_prob} admits a $\ell_{x^\star}\ell_{\tilde{\xi}}$-Lipschitz solution map $\tilde{x}^\star$, in $\psi$. This holds since the composition of Lipschitz maps is also Lipschitz.

Based on the special problem classes discussed in the previous subsection, we can also derive the following result which incorporates non-differentiable perturbations with tighter bound on the Lipschitz constant than derived from \eqref{eq:time_varying_prob}:
\begin{corollary}\label{cor:special_case}
    Consider the problem
    \begin{align}\label{eq:prob_spec}
    \begin{split}
        \minimize \quad &\hat{f}(x - c(\xi)) \\
        \st \quad &  U x \leq v(\xi) \,
    \end{split}
    \end{align}
    where $\xi \in  \Xi \subset \R^r$. Define $\hat{\calX}(\xi) := \{ x \, | \, U x \leq v(\xi) \}$ and let
    \begin{enumerate}[label=(\roman*)]
        \item  $\hat{f}: \R^n \rightarrow \R$ be twice continuously differentiable, $\alpha$-strongly convex, with $\beta$-Lipschitz $\nabla \hat{f}$,
        \item $c: \Xi \rightarrow \R^n$ and $v: \Xi \rightarrow \R^m$ be $\ell_c$- and $\ell_v$-Lipschitz, respectively, and
        \item  $U \in \R^{n \times m}$ and $\omega > 0$ such that for every $\xi \in \Xi$ and ever $x \in \hat{\calX}(\xi)$ one has $\omega^2 \bbI \preceq U_{\bfI_{x}^\xi} U_{\bfI_{x}^\xi}^\T$.
    \end{enumerate}
If $\hat{\calX}(t) \neq \emptyset$ for every $\xi \in \Xi$, then the primal and dual solution maps $x^\star: \Xi \rightarrow \R^n$ and $\mu^\star: \Xi \rightarrow \R^m$ of \eqref{eq:prob_spec} are Lipschitz with respective bounds on Lipschitz constants
\begin{align*}
    \ell_{x^\star} =  \sqrt{\frac{\beta}{\alpha}} \left(  \frac{\beta \ell_c}{\alpha} +\frac{ \ell_v}{ \omega} \right) \, \, \text{and} \, \,
    \ell_{\mu^\star} =\frac{\beta^{3/2}}{\alpha^{1/2} \omega} \left(  \frac{\beta \ell_c}{\alpha} +\frac{ \ell_v}{ \omega} \right) \, .
\end{align*}
\end{corollary}
\begin{proof}
    The key point about \cref{cor:special_case} is the fact that $c$ and $v$ are Lipschitz continuous, but not twice differentiable. To address this issue, consider,  instead of \eqref{eq:prob_spec}, the problem
    \begin{align}\label{eq:prob_spec_red}
    \begin{split}
        \minimize \quad &\hat{f}(x - c) \qquad
        \st \quad   U x \leq v \,
    \end{split}
    \end{align}
    parametrized in $(c, v)$. Under the stated assumptions, the parametrized objective $(x, c) \mapsto \hat{f}(x- c)$ and constraint function $(x, v) \mapsto Ux - v$ are twice continuously differentiable in $(x,c)$ and $(x,v)$, respectively. Therefore, the solution map $(c, v) \mapsto \hat{x}^\star(c,v)$ of \eqref{eq:prob_spec_red} is Lipschitz by \cref{cor:ell.lipschitz.optimizer}. In fact, by considering the reasoning for translational perturbations of the objective and righthand perturbations of polyhedra in the previous subsection, we know that $\hat{x}^\star$ is Lipschitz in $c$ and $v$ with Lipschitz constants $(\tfrac{\beta}{\alpha})^{3/2}$ and  $\sqrt{\tfrac{\beta}{\alpha}}\tfrac{1}{\omega}$, respectively.

    Moreover, replacing $c$ with $c(\xi)$ and using $\ell_c$-Lipschitz continuity of $c$ we have, for all $c',c \in \R^n$ and all $v \in \R^m$,
    \begin{align}\label{eq:c_lip}
        \norm{ \hat{x}^\star(c(\xi'), v) -  \hat{x}^\star(c(\xi), v)} \leq \sqrt{\tfrac{\beta}{\alpha}} \tfrac{\beta}{\alpha} \ell_c
        \norm{\xi' - \xi}
    \end{align}
    Similarly, replacing $v$ with $v(\xi)$ we have
    \begin{align}\label{eq:v_lip}
        \norm{  \hat{x}^\star(c, v(\xi')) -  \hat{x}^\star(c, v(\xi))} \leq \sqrt{\tfrac{\beta}{\alpha}} \tfrac{1}{\omega} \ell_v \norm{\xi' - \xi}
    \end{align}
    for all $c \in \R^n$ and all $v', v \in \R^m$.

    The mutual uniformity of the these Lipschitz constants, i.e., the fact that \eqref{eq:c_lip} and \eqref{eq:v_lip} are independent of $v$ and $c$, respectively, allows us to write
    \begin{align*}
         \norm{ x^\star(\xi') -  x^\star(\xi)}    & \leq \sqrt{\tfrac{\beta}{\alpha}} \left(  \tfrac{\beta \ell_c}{\alpha} +\tfrac{ \ell_v}{ \omega} \right) \| \xi' - \xi \| \, ,
    \end{align*}
    where $x^\star(\xi) = \hat{x}^\star(c(\xi), v(\xi))$ and thus the solution map $x^\star: \Xi \rightarrow \R^n$ of \eqref{eq:prob_spec} is Lipschitz with the desired Lipschitz constant. The Lipschitz constant for $\mu^\star$ follows analogously.
\end{proof}
\section{Tracking Performance of Continuous-time Algorithms for Time-Varying Optimization}\label{sec:tracking.performance}
In this section, we derive a tracking result that applies to \emph{any} continuous-time algorithm with absolutely continuous trajectories that satisfy a common monotonicity-like property.
Although the same property has previously been exploited for specific algorithms, our tracking result is more general since it does not require any additional assumptions on the (running) algorithm. For instance, our performance bound can be applied to algorithms whose trajectories are Krasovskii \cite[Ch. 5]{goebelHybridDynamicalSystems2012} or Filippov \cite{filippovDifferentialEquationsDiscontinuous1988} solutions of differential inclusions.
In the forthcoming \cref{sec:sweeping.gradient.descent}, we require the generality of this result for our novel discontinuous sweeping gradient flow.

Generally, in this section, we consider a time-varying optimization problem with non-smooth variations in time whose instantaneous optimizer the running algorithm tracks. Hence, instead of periodically referring back to the sensitivity results of the previous sections, we make the following assumption:
\begin{assumption}\label{ass:tracking_opt}
    The problem \eqref{eq:basic_opt_prob} with $\xi=\txtold$ and $\tilde{\xi}: \R_{\geq 0} \rightarrow \R^r$ admits a unique (global) optimizer for every $t \in \R_{\geq 0}$ and the solution map $\txs:\R_{\geq 0} \rightarrow \R^n$ is $\ell$-Lipschitz.
\end{assumption}

The following proposition guarantees the bounded asymptotic behavior of the running algorithm, i.e., the distance $\norm{x(t)-\txs(t)}$ between the trajectories of the running algorithm $x(t)$ and instantaneous optimizer $\txs(t)$ is asymptotically bounded.

\begin{proposition}\label{cor:tracking.performances}
    Consider  \eqref{eq:basic_opt_prob} with $\xi=\txtold$ and $\tilde{\xi}: \R_{\geq 0} \rightarrow \R^r$ and let \cref{ass:tracking_opt} be satisfied. Furthermore, let $x: \R_{\geq 0} \rightarrow \R^n$ be absolutely continuous and assume that
  \begin{align}\label{eq:monot_cond_grad}
       \left\langle \dot x(t), x(t) - \txs(t) \right\rangle  \leq - a \| x(t) - \txs(t) \|^2,
\end{align}
    holds for almost every $t \in \R_{\geq 0}$ where $a \in \R_{>0}$. Then, $\underset{t \rightarrow \infty} \lim \sup \| x(t) - \txs(t) \| \leq \tfrac{\ell}{a}$ holds.
    Moreover, if $\| x(0) - \txs(0) \|\leq \tfrac{\ell}{a}$, then $\|x (t) - \txs(t) \| \leq \tfrac{\ell}{a}$ holds for all $t \geq 0$.
\end{proposition}

In the next subsection we provide the necessary technical results and the proof of \cref{cor:tracking.performances}.

\subsection{Proof of \cref{cor:tracking.performances}}
Before we prove the proposition, we derive our key technical result from Barbalat's lemma which is at the base for invariance-like theorems \cite{khalilNonlinearSystems2002,fischerLaSalleYoshizawaCorollariesNonsmooth2013}.

Recall that the function $f:S \rightarrow \R$ is said to be \emph{uniformly continuous on $ S \subseteq \R$} if and only if for all $\epsilon>0$ there exists $\delta>0$ such that for all $ x', x \in S$
\begin{align*}
   \abs{x'-x}<\delta \implies \abs{f(x')-f(x')}<\epsilon.
\end{align*}

Hence, Barbalat's lemma reads as follows:

\begin{lemma}[Barbalat's lemma] \cite[Lemma 8.2]{khalilNonlinearSystems2002}\label{lem:barbalat} Let $W:\R_{\geq 0} \rightarrow \R$ be uniformly continuous. If ${\lim}_{t\rightarrow \infty} \int_{0}^t W(\tau) d\tau$ exists and is finite, then $\lim_{t \rightarrow \infty} W(t) = 0$.
\end{lemma}
Moreover, we also require the following lemma:
\begin{lemma}\label{lem:unif_monot}
    If $W:\R_{\geq 0} \rightarrow \R$ is continuous, lower bounded, and non-increasing, then $W$ is uniformly continuous.
\end{lemma}

\begin{proof}
    Let $L = \inf_{t \geq 0} W(t)$ and consider any $\epsilon > 0$. Let $T > 0$ be such that $W(t) \leq L + \epsilon$ for all $t > T$.  Since a continuous function is absolutely continuous on a compact interval, $W$ is uniformly continuous on the compact interval $[0, T']$ with $T' > T$. In particular, for the given $\epsilon$, there exists $\delta > 0$ such that $|t - \tau | \leq \delta$ implies $| W(t) - W(\tau) | \leq \epsilon$ for any $t,\tau \in [0,T']$. Without loss of generality, let $\delta < T' - T$.
    Next, note that for any $t, \tau \in [T, \infty)$ (and, in particular, for $| t - \tau | \leq \delta$), we have that $| W(t) - W(\tau) | \leq \epsilon$.
    Therefore, for any $t, \tau \geq 0$ with $| t - \tau | \leq \delta$ we have $| W(t) - W(\tau) | \leq \epsilon$. Since $\epsilon$ is arbitrary, $W$ is uniformly continuous.
\end{proof}
Finally, we can present the key technical result necessary for the proof of \cref{cor:tracking.performances}:
\begin{theorem} \label{lem:tracking.performances}
  Let $y: \R_{\geq 0} \rightarrow \R^n$ be absolutely continuous such that for almost every $t \geq 0$, some $a > 0$, and $b \geq 0$ it holds
    \begin{align}
        \left\langle \dot y(t), y(t) \right\rangle  \leq - a \norm{ y(t)}^2  + b \norm{ y(t) }, \label{eq:bounded.inner.product}
    \end{align}
    Then, $\underset{t \rightarrow \infty}{\lim \sup} \, \norm{y(t)}  \leq \tfrac{b}{a}$ holds.
    Moreover, if $\norm{y(T)}\leq \frac{b}{a}$ holds for some $T \geq 0$, then $\norm{y(t)}\leq \frac{b}{a}$ holds for all $t\geq T$.
\end{theorem}
\begin{proof}
First we note that the absolutely continuous function on a compact interval has a bounded variation and hence, it is differentiable almost everywhere. Moreover, we define $V:\R_{\geq 0} \rightarrow \R_{\geq 0}$ as
    \begin{align*}
       V(t)\coloneqq \max \left\{  \tfrac{1}{2}\left(\norm{y(t)}^2 - (\tfrac{b}{a})^2 \right) , 0 \right\} \, ,
    \end{align*}
and let $W:\R_{\geq 0} \rightarrow \R_{\geq 0}$ be given as
\begin{align*}
    W(t) \coloneqq \max \left\{ a\norm{y(t)} \left( \norm{y(t)}- \tfrac{b}{a} \right), 0 \right\}\, .
\end{align*}
Note that both $V$ and $W$ are absolutely continuous. Furthermore,  $V(t) > 0$ and $W(t) > 0$ holds for all $t$ for which $\norm{y(t)} > \tfrac{b}{a}$ and $V(t) = W(t) = 0$ otherwise.
Moreover, for almost all $t \geq 0$ for which $\norm{y(t)} \geq \tfrac{b}{a}$ we have
\begin{align*}
    \begin{split}
   \dot{V}(t) & = \left \langle \dot{y}(t), y(t) \right \rangle
    \leq  - a \norm{y(t)}^2 + b \norm{y(t)} = - W(t) \leq 0 \,.
    \end{split}
\end{align*}
Otherwise, for almost all $t$ for which $\norm{y(t)} < \tfrac{b}{a}$, we have $\dot{V}(t)= W(t) = 0$. Hence, the growth of $V$ is bounded by $\dot{V}(t)\leq -W(t)\leq 0$ for almost all $t \geq 0$.

To prove the lemma, we apply \cref{lem:barbalat} to $W$. Hence, we need to show that \textit{(i)} the integral $\underset{t\rightarrow \infty}{\lim} \int_{0}^t W(\tau) d\tau$ exists and is finite, and that \textit{(ii)} $W$ is uniformly continuous.

To show \textit{(i)}, we exploit the fact that the growth of $V$ is bounded by $-W$ to establish the bound
\begin{align*}
    \int_{0}^t W(\tau) d \tau \leq - \int_{0}^t \dot{V}(\tau) d \tau  = V(0)- V(t)
    \leq V(0) \, , 
\end{align*}
where the last inequality follows from the fact that $V(t)$ is non-negative for all $t \in \R_{\geq 0}$, by definition.
Further, the left-hand side integral is non-decreasing in $t$ (since $W(t) \geq 0$ for all $t$). Hence, as $t \rightarrow \infty$, the limit exists and is finite.

For \textit{(ii)} it suffices to show that $W$ is non-increasing itself. Then, since $W$ is continuous and lower bounded, it follows from \cref{lem:unif_monot} that $W$ is uniformly continuous.

For this purpose, simply note that, for almost all $t$ for which $\norm{y(t)} \leq b/a$, we have $\dot{W}(t)=0$ and, for almost all $t$ for which $\norm{y(t)} > b/a$, it holds that
\begin{align*}
\dot{W}(t) &= \underbrace{\left(2 a -\tfrac{ b }{\norm{y(t)}} \right)}_{\geq a > 0} \underbrace{\langle \dot{y}(t),y(t) \rangle}_{= \dot V(t) < 0} < 0 \, .
\end{align*}

Thus, we have shown that $W$ is uniformly continuous and $\underset{t\rightarrow \infty}{\lim} \int_{0}^t W(\tau) d\tau$ exists. Therefore, from \cref{lem:barbalat}, it follows that $W(t) \rightarrow 0$ for $t \rightarrow \infty$ and consequently, from the definition of $W$, we have $\lim \underset{t \rightarrow \infty}{\sup} \| y (t) \| \leq b / a$.

We show the last part of the lemma by contradiction. Assume that there exists $T\geq 0$ such that $\norm{y(T)}\leq \frac{b}{a}$ and that there exists $t_1\geq T$ such that $\norm{y(t_1)}>\frac{b}{a}$. That means $V(t_1)>V(T)=0$, but $\dot{V}(T)=0$, hence $V(T+\delta)=0$, for all $\delta \geq 0$. Consequently $V(t_1)=0$ which is contradiction to $V(t_1)>V(T)$, hence $\norm{y(t)}\leq \frac{b}{a}$, for all $t\geq T$.
\end{proof}

Note that the second part of \cref{lem:tracking.performances}, that guarantees the invariance of $\tfrac{b}{a} \bbB$, is essentially a consequence of Nagumo's theorem \cite[Thm. 1.2.1]{aubinViabilityTheory1991}. Likewise, \cref{lem:tracking.performances} and its proof relate to Lyapunov functions, barrier functions, etc.

Lastly, in order to show that \cref{cor:tracking.performances} follows from \cref{lem:tracking.performances}, we  first note  $x$ and $\txs$ are both (individually) differentiable for almost all $t \geq 0$ and thus they are (jointly) differentiable for almost all $t \geq 0$ (since the union of two zero-measure sets has measure zero).
Hence, the instantaneous optimizer trajectory $\txs$ satisfies \cref{ass:tracking_opt} for almost all $t \geq 0$, and therefore, by Cauchy-Schwarz, it holds that $ |\langle \dot{\tilde{x}}^\star(t),x(t)-\txs(t)\rangle| \leq \ell \norm{x(t)-\txs(t)}$ for almost all $t \geq 0$. Using \eqref{eq:monot_cond_grad} and letting $y(t) \coloneqq x(t)-\txs(t)$ (which is absolutely continuous) and $b\coloneqq \ell$, the inequality \eqref{eq:bounded.inner.product} is satisfied and the proof of \cref{cor:tracking.performances} follows immediately.

The \emph{strong monotonicity} property \eqref{eq:monot_cond_grad} is, for example, satisfied with $a = \alpha$ for an unconstrained gradient flow $\dot x = - \nabla_x f(x,t)$ if $f$ is $\alpha$-strongly convex for all $t \geq 0$, which is further detailed in \cref{sec:tv.u.o}.
\section{Continuous-Time Constrained Optimization} \label{sec:sweeping.gradient.descent}
Before presenting a novel continuous-time optimization scheme, we review some results regarding \emph{sweeping processes}.

\subsection{Perturbed Sweeping Processes}

 Given a non-empty convex set $\calX$, recall that the normal cone of $\calX$ at $x\in \calX$ is given by
    \begin{align} \label{eq:normal.cone}
        N_x \calX := \cl \{ \eta \, | \, \forall y \in \calX: \, \left\langle \eta, y - x \right\rangle \leq 0 \} \, .
    \end{align}
Sweeping processes \cite{moreauEvolutionProblemAssociated1977,thibaultSweepingProcessRegular2003,kunzeIntroductionMoreauSweeping2000} have originally been formulated to describe a sweeping effect of a moving impenetrable boundary on a mobile object. When the object itself is subject to a perturbation (e.g. a force), we use a \emph{perturbed sweeping process} \cite{castaingEvolutionProblemsAssociated1996, edmondRelaxationOptimalControl2005} which is formally defined as
\begin{align}\label{eq:sweeping.process}
    \dot{x}(t) \in \Phi(x,t) - N_x \calX(t), \ \text{with} \ x(0) \in \calX(0) \, .
\end{align}
where $\Phi: \R^n \times \R \rightarrow \R^n$ is a time-varying vector field, and, for every $t \geq 0$, $\calX(t) \subset \R^n$ is a non-empty closed, convex set.
The idea behind \cref{eq:sweeping.process} is that $x$ evolves according to $\Phi$ in the interior of $\calX(t)$ (where $N_x \calX(t) = \{ 0 \}$), but when $x$ touches the boundary, a normal ``force'' is exerted on $x$ to keep the state feasible, even as set $\calX(t)$ varies.

A \emph{solution of \eqref{eq:sweeping.process}} is an absolutely continuous map $x: [0, T] \rightarrow \R^n$ for some $T > 0$ such that, for almost all $t\in [0, T]$ \eqref{eq:sweeping.process} is satisfied, and $x(\txt) \in \tilde{\calX}(\txt)$ holds for all $t \in [0,T]$. A \emph{complete solution} of \eqref{eq:sweeping.process} in an absolutely continuous map $x: \R_{\geq 0} \rightarrow \R^n$ such that the restriction to any compact interval $[0, T]$ for $T > 0$ is a solution of \eqref{eq:sweeping.process}.

The following theorem is based on \cite[Thm. 1]{edmondRelaxationOptimalControl2005} and simplified\footnote{%
The original result applies to prox-regular sets in a Hilbert space instead of convex sets in $\R^n$. Furthermore, in the original theorem \cref{hypothesis:H2} is more general and requires $\calX$ to vary in a \emph{absolutely continuous} way instead of a \emph{Lipschitz continuous} way. Finally, in the original statement \cref{hypothesis:H3} is split into a \emph{local Lipschitz} and \emph{linear growth} condition. For simplicity, we assume global Lipschitz continuity to meet both requirements.} for our specific purposes.

\begin{theorem} \label{thm:exist_uniq_sweeping_process}
For some $T>0$, let $\calX : [0, T] \rightrightarrows \R^n$ satisfy
\begin{romcases}
    \item \label{hypothesis:H1} $\calX(t)$ is nonempty, closed and convex for all $t \in [0, T]$,
    \item \label{hypothesis:H2} $\ell_\calX > 0$ exists such that for any $z \in \R^n$ and $t', t \in [0, T]$,
    \begin{align}\label{eq:sweep_set_variation}
        |d(z,\calX(t'))-d(z, \calX (t))| \leq  \ell_\calX |t'- t|,
    \end{align}
    where $d$ denotes the point-to-set distance, and it is given by $d(x,\calX(t)) \coloneqq \inf \{\norm{x-y} \, : \, y \in \calX(t) \}$.
    \end{romcases}
    Further, let $\Phi: [0, T] \times \R^n \rightarrow \R^n$ be measurable in $t$ and
    \begin{romcases}[resume]
        \item \label{hypothesis:H3} there exists an integrable function $\gamma: [0, T] \rightarrow \R_{\geq 0}$ such that $\norm{\Phi(x,t)-\Phi(y,t)}\leq \gamma(t) \norm{x-y} $ holds for all $t \in [0, T]$ and any $x, y \in \R^n$.
    \end{romcases}
    Then, for any initial condition $x(0) \in \calX(0)$, the \emph{perturbed sweeping process}
     \begin{align*}
         \dot{x} \in \Phi(x,t )-N_{x} \calX(t) \, \quad x \in \calX(t) .
     \end{align*}
     admits a unique solution $x: [0, T] \rightarrow \R^n$.
\end{theorem}
\subsection{Sweeping Gradient Flows}
As a new running algorithm for constrained time-varying optimization we consider a perturbed sweeping process of the form \eqref{eq:sweeping.process} where $\Phi(x,t) \coloneqq -\grad_x f(x,t)$ is the negative gradient of a time-varying cost function $f: \R^n \times \R \rightarrow \R$, thus the resulting \emph{sweeping gradient flow} is given by
\begin{align}\label{eq:sweeping.gradiend.descent}
  \dot{x} (t) \in  - \grad_x {f}(x,\txt) - N_{x} {\calX}(\txt) \, , \quad \ x(0
  )\in {\calX}(0).
\end{align}

In the following, we consider sweeping gradient flows for the special problem structure from \cref{cor:special_case} for which we can give easy-to-interpret results. However, \eqref{eq:sweeping.gradiend.descent} is well-defined for much more general setups.\footnote{For instance, existence of trajectories for sweeping gradient flows is guaranteed even if $\calX$ is not convex, but \emph{prox-regular}, as long as $\nabla_x f$ is measurable in $t$, $\calX(t)$ is non-empty for all $t \geq 0$ and, roughly-speaking, $\calX$ varies in an absolutely continuous way~\cite{edmondRelaxationOptimalControl2005}.
To guarantee a global asymptotic bound on the tracking error, however, we generally require \cref{ass:feasible,ass:cvx,ass:unif_licq,ass:dual_bound} to hold such that \cref{cor:ell.lipschitz.optimizer,cor:tracking.performances} can be applied.
}

The key insight from this section is that our sensitivity bounds and the generalized tracking results of the previous sections can be applied to non-trivial discontinuous optimization dynamics involving time-varying constraints. Namely, we establish the following result:

\begin{theorem}\label{thm:special_case_tracking}
Consider the problem
\begin{align}\label{eq:spec_case_prob_tracking}
\begin{split}
    \minimize \quad &\hat{f}(x - c(t)) \\
    \st \quad &  U x \leq v(t) \, .
\end{split}
\end{align}
where $t \geq 0$. Define $\hat{\calX}(t):= \{ x\, | \,U x \leq v(t) \}$ and let
    \begin{enumerate}[label=(\roman*)]
        \item  $\hat{f}: \R^n \rightarrow \R$ be twice continuously differentiable, $\alpha$-strongly convex, and $\nabla \hat{f}$ is $\beta$-Lipschitz,
        \item $c: \R_{\geq 0} \rightarrow \R^n$ and $v: \R_{\geq 0} \rightarrow \R^m$ be $\ell_c$- and $\ell_v$-Lipschitz continuous,
        \item  $U \in \R^{n \times m}$ and $\omega > 0$ such that for every $\xi \in \Xi$ and ever $x \in \hat{\calX}(\xi)$ one has $\omega^2 \bbI \preceq U_{\bfI_{x}^\xi} U_{\bfI_{x}^\xi}^\T$.
    \end{enumerate}
If $\hat{\calX}(t) :=  \{ x \, | \, U x \leq v(t) \} \neq \emptyset$ for all $t \geq 0$, then the sweeping gradient flow
\begin{align} \label{eq:sweep_grad_spec}
    \dot x \in - \grad \hat{f}(x - c(t)) - N_x \hat{\calX}(t)
\end{align}
admits a complete solution $x: \R_{\geq 0} \rightarrow \R^n$  for every initial condition $x(0) \in \hat{\calX}(0)$ and it holds that
\begin{align} \label{eq:thm.special.tracking.bound}
    \underset{t \rightarrow \infty}{\lim \sup} \| x(t) - x^\star(t) \| \leq \frac{\ell_t}{\alpha} :=\frac{\beta^{1/2}}{\alpha^{3/2}} \left(  \frac{\beta\ell_c}{\alpha}  +\frac{ \ell_v}{ \omega} \right) \, ,
\end{align}
where $x^\star(t)$ is the unique optimizer of \eqref{eq:spec_case_prob_tracking} at time $t$.

Furthermore, if $\| x(0) - x^\star(0) \| \leq \tfrac{\ell_t}{\alpha}$, then $\| x(t) - x^\star(t) \| \leq \tfrac{\ell_t}{\alpha}$ holds for all $t \geq 0$.
\end{theorem}

To show \cref{thm:special_case_tracking} we first need to prove the existence of complete solutions for \eqref{eq:sweeping.gradiend.descent}.

\begin{lemma} \label{prop:existence_sweep_grad}
Given the setup of \cref{thm:special_case_tracking}, \eqref{eq:sweep_grad_spec} admits a unique complete solution for any initial point $x(0) \in \hat{\calX} (0)$.
\end{lemma}

Even though \cref{prop:existence_sweep_grad} is an existence result, we can prove it using the sensitivity results of the previous section.
\begin{proof}
The proof is an application of \cref{thm:exist_uniq_sweeping_process}. In order to apply \cref{thm:exist_uniq_sweeping_process} to \eqref{eq:sweep_grad_spec}, note that for the setup in \cref{thm:special_case_tracking} the requirements \cref{hypothesis:H1} in \cref{thm:exist_uniq_sweeping_process} (non-empty closed convex $\hat{\calX}(t)$) and \cref{hypothesis:H3} (Lipschitz vector field) hold by assumption on any compact interval $[0, T] \subset \R_{\geq 0}$. Also, $ \grad f(x,t) =  \grad \hat{f}(x - c(t))$ is measurable in $t$ since it is continuous in $t$.
It remains to show \cref{hypothesis:H2} by establishing \cref{eq:sweep_set_variation}.

For this purpose, we consider $d(z, \calX(t))$ as the solution of the parametrized optimization problem
\begin{align} \label{eq:min.distance.optimization.problem}
\begin{split}
    \underset{y}{\minimize} \quad &  \tfrac{1}{2} \norm{y-z}_2^2,\\
    \st \quad &  U y \leq v(t) \, ,
    \end{split}
\end{align}
where $t$ varies but $z$ is fixed. As a problem parametrized in $t$, \cref{eq:min.distance.optimization.problem} falls into the class of problems to which \cref{cor:special_case} applies. Namely, we have $\alpha = \beta = 1$ (since $\hat{f}(x) := \tfrac{1}{2}\| x \|^2$), $\ell_c = 0$ (since $c \equiv 0$), and $\omega$ as defined in \cref{thm:special_case_tracking}. Further, by the assumption in \cref{thm:special_case_tracking}, the feasible set of \eqref{eq:min.distance.optimization.problem} (which is equivalent to $\hat{\calX}(t)$) is non-empty for all $t \geq 0$.

Therefore, by \cref{cor:special_case}, the solution map $t \mapsto y^\star_z(t)$ of \cref{eq:min.distance.optimization.problem} is Lipschitz continuous with a bound on the Lipschitz constant that is independent of $z$ and given by $\ell_{y^\star} = \ell_v / \omega$. Therefore, we have $| y^\star_z(t') - y^\star_z(t) | \leq \ell_{y^\star} | t' - t |$ for any $z \in \R^n$ and $t, t' \geq 0$, and \cref{hypothesis:H2} in \cref{thm:exist_uniq_sweeping_process} holds on any compact interval $[0, T]$.

Hence, \cref{thm:exist_uniq_sweeping_process} guarantees the existence of a unique solution of \eqref{eq:sweep_grad_spec} for every initial condition $x(0) \in \hat{\calX}(0)$, and for every compact interval $[0,T]$ and hence, by definition, a complete solution on $\R_{\geq 0}$.
\end{proof}

\begin{proof}[Proof of \cref{thm:special_case_tracking}]
By \cref{cor:special_case}, \eqref{eq:spec_case_prob_tracking} admits a unique primal solution $x^\star(t)$ for every $t \geq 0$.
In particular, the solution map $t \mapsto x^\star(t)$  of \eqref{eq:spec_case_prob_tracking} is $\ell_t$-Lipschitz with $\ell_t :=   \sqrt{\frac{\beta}{\alpha}} \left(  \frac{\beta\ell_c}{\alpha}  +\frac{\ell_v}{ \omega}  \right)$. Consequently, \cref{ass:tracking_opt} is satisfied.
Furthermore, \cref{prop:existence_sweep_grad} guarantees the existence of a unique complete solution $x : \R_{\geq 0} \rightarrow \R^n$ of \eqref{eq:sweep_grad_spec} for every initial condition $x(0) \in \hat{\calX}(0)$.
Next, we verify that \eqref{eq:monot_cond_grad} holds, i.e., we show that $\left \langle \dot{x}(t), x(t) - x^\star(t) \right\rangle \leq -\alpha \| x(t) - x^\star(t) \|^2$.

In the following let $f(x, t) := \hat{f}(x - c(t))$.
Recall that, for almost all $t \geq 0$, we have $\dot{x}(t)=-\grad_x  f(x(t), t) - \eta(t)$ for some $\eta(t) \in N_{x(t)} \hat {\calX}(t)$. Further, because $x^\star(t)$ solves \eqref{eq:spec_case_prob_tracking}, it satisfies the KKT conditions which is equivalent to saying that $-\grad_x f(x^\star(t) ,t) - \eta^\star(t)=0$ for some $\eta^\star(t) \in N_{x^\star(t)} \hat{\calX}(t)$. Putting these two insights together, we have that
\begin{equation*}
\dot{x}(t)=-\grad_x  f(x(t),t) - \eta(t) + \grad_x  f(x^\star(t),t) + \eta^\star(t) \,.
\end{equation*}

Next, we omit the argument $t$ for $x, \ \txs$, $\eta, \ \eta^\star$ and $\hat{\calX}$. By definition \eqref{eq:normal.cone} of the normal vectors $\eta \in N_{x} \hat{\calX}$ and $\eta^\star \in N_{x^\star} \hat{\calX}$, for $\hat{\calX}$, we have $\langle \eta ,x-x^\star \rangle \geq 0$ and $\left \langle \eta^\star, x^\star - x \right\rangle \geq 0$, respectively. Since $f$ is $\alpha$-strongly convex we have
\begin{align}\label{eq:strong.convexity.property2}
        & \left\langle \grad f(x,t)-\grad f(x^\star,t), (x-x^\star) \right\rangle \geq \alpha\norm{x- x^\star}^2 \,.
\end{align}
Combining these facts, we get
\begin{align*}
     \left\langle \dot{x}, x\! -\! x^\star \right\rangle &\! =\! \left\langle- (\grad \! f(x,t) \!- \! \grad\! f(x^\star,t)) \!-\! ( \eta\!-\! \eta^\star), x\!-\! x^\star \right\rangle \\
     &\!=\! -\! \left\langle\eta, x\!-\! x^\star \right\rangle \!-\! \left\langle \eta^\star, x^\star\!-\! x \right\rangle  \\
      & \ \ - \! \left\langle\grad \! f(x,t) \!- \! \grad\! f(x^\star,t), x\!-\! x^\star \right\rangle\\
    & \!\leq\! - \alpha \norm{x - x^\star}^2  .
\end{align*}

Consequently, by taking $a=\alpha$, \eqref{eq:monot_cond_grad} holds. Thus, \cref{cor:tracking.performances} is applicable and yields the desired asymptotic tracking bound and completes the proof.
\end{proof}

\section{Illustrative, Numerical Examples} \label{sec:numerical.examples}
In this section, we provide two simple numerical examples to illustrate our results for simple time-varying optimization setups. First, we consider an unconstrained optimization problem in one dimension, with non-smooth change in time to demonstrate the tightness of our bound in that case. Second, we consider a constrained time-varying optimization problem to illustrate the behavior of sweeping gradient flows.

\subsection{Time-Varying Unconstrained Optimization} \label{sec:tv.u.o}

Consider the problem of minimizing
\begin{align*}
    f(x,\tilde{\xi}(t))&= \norm{x-\tilde{\xi}(t)}^2,
\end{align*}
where $\tilde{\xi}$ is a triangular wave (hence non-smooth) with period $\tau=4$s. That is, for all $n\in \N$, $\tilde{\xi}$ is defined as
   \begin{align*}
   \tilde{\xi}(t)&=\left \{ \begin{array}{cc}
       \!(t-\tau(n-1))-1,  & \!t \in [\tau(n-1),\tau(n-0.5) \\
       \!-(t-\tau(n-1)) +3,  & \! t \in [
       \tau(n-0.5),\tau n)
    \end{array} \!\!  \right. .
    \end{align*}
In particular, in reference to \eqref{eq:time_varying_prob} we have $\psi=t$ and $f(x, \xi) := \norm{x - \xi}^2$ is $\alpha$-strongly convex  with $\alpha=2$. Furthermore, the solution map given by $\tilde{x}^\star(t)=\txtold$ is $\ell_{\tilde{\xi}}$-Lipschitz with $\ell_{\tilde{\xi}}=1$.
From \eqref{eq:unconstrained.ell.constant} we have $\sup_{\xi \in [-1,1]}\norm{\nabla^2_{\tilde{\xi} x}f}= 2$, the solution trajectory $x^\star$ is $\ell_{x^\star}$-Lipschitz with respect to $\xi$, where $\ell_{x^\star}=1$. Therefore, the solution trajectory $\txs$ is  $\ell_t$-Lipschitz continuous in $t$ with $\ell_t := \ell_{x^\star}\ell_{\tilde{\xi}}=1$.

For all $t\geq 0$, the unconstrained gradient flow is given by
\begin{align*}
    \dot{x}(t)=-\grad_x f(x(t),\txtold)=-2(x(t)-\txtold).
\end{align*}

Finally, combining the zero cost function gradient, at the unconstrained optimizer (i.e.,  $\grad_x f(\txs(t),\txtold)=0$) with the $\alpha$-strongly convex property \eqref{eq:strong.convexity.property2}, for almost all $t\geq 0$ we have
\begin{align*}
    \langle \dot{x}, x(t)-\txs(t) \rangle \leq -\alpha \norm{ x(t)-\txs(t)}^2.
\end{align*}
Hence, from \cref{cor:tracking.performances}, ${\lim \sup}_{t \rightarrow \infty} \norm{x(t)-\txs(t)} \leq 0.5$.

\begin{figure}[tb]
\centering
\includegraphics[width=.98\columnwidth]{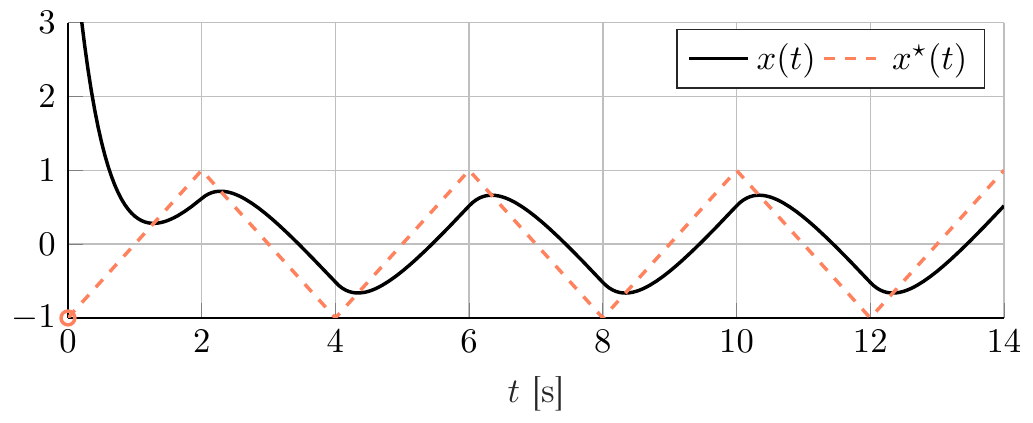}
\vspace{-1em}
\caption{Trajectories of the gradient descent $x(t)$ (black, solid) and the unconstrained optimizer ${x}^\star(t)$ (vermilion, dashed). \label{fig:trajectories1d}}
\vspace{1em}
\includegraphics[width=.98\columnwidth]{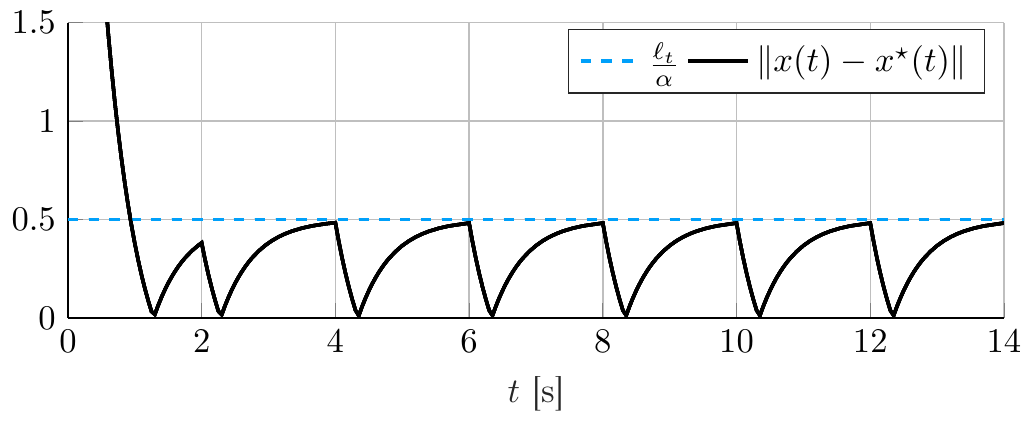}
\vspace{-1em}
\caption{Tracking error $\norm{x(t)-{x}^\star(t)}$ (black, solid), and the tracking error upper bound $\tfrac{\ell_t}{\alpha}$ (blue, dashed). \label{fig:trajectories1d.tracking.error}}
\end{figure}

\cref{fig:trajectories1d} shows the trajectory $x(t)$ (black, solid) obtained from the gradient descent algorithm with initial condition $x(0)=5$ and the instantaneous optimizer trajectory $\txs(t)$ (vermilion, dashed).

\cref{fig:trajectories1d.tracking.error} shows the distance between the trajectories $\norm{x(t)-\txs(t)}$ (black, solid), and the analytical sensitivity bound $\tfrac{\ell_t}{\alpha}=0.5$ (blue, dashed). Easily, it can be noticed that as $t\rightarrow \infty$, the distance always stays within the analytical sensitivity bound. Additionally, in this example, the analytical sensitivity bound is indeed tight. Moreover, this example perfectly illustrates the invariance property of a ball of radius  $\tfrac{\ell_t}{\alpha}$ around the time-varying unconstrained optimizer $\txs(t)$. Namely, the algorithm trajectory never perfectly tracks the unconstrained optimizer trajectory, yet always stays within $\tfrac{\ell_t}{\alpha}$ from the optimizer.
\subsection{Time-Varying Constrained Optimization}

In this example, we consider the same problem setup as in \cref{thm:special_case_tracking} with smooth changes in time. Namely, let
\begin{align*}
    Q\coloneqq\begin{bmatrix}
    12 & -8 \\-8 & 10
    \end{bmatrix} \succ 0, \ P(t)\coloneqq-\begin{bmatrix}
    3\sin(t+3)+1.3t \\ 2\tanh(t-3)+0.71t
    \end{bmatrix}, \
\end{align*}
let $c(t)\coloneqq - \tfrac{1}{2}Q^{-1} P(t)$, and the cost function is defined as
\begin{align*}
    \hat{f}(x(t)-c(t))&= (x(t) + \tfrac{1}{2}Q^{-1} P(t))^\T Q (x(t) + \tfrac{1}{2}Q^{-1} P(t)).
\end{align*}
In particular, $\hat{f}$ is $\alpha$-strongly convex with $\alpha=2\sigma^{\min}(Q)$, and gradient $\grad\hat{f}$ is $\beta$-Lipschitz with $\beta=2\sigma^{\max}(Q)$. Furthermore, we can compute $\ell_c =\sup_{t \geq 0}\norm{\dot{P}(t)}$.

Next, the constraint set is given by
\begin{align*}
    \hat{\calX}(t)=\{ x(t) \in \R^n \ | \ Ux(t)\leq v(t) \},
\end{align*}
where $v(t)\coloneqq (V_1 t +V_2)$ and
\begin{align*}
    U\coloneqq \begin{bmatrix}
   -2&1\\1&-1\\0.5&1\\-3&-1
    \end{bmatrix}, \quad V_1 \coloneqq \begin{bmatrix}
    -0.05\\-0.3\\0.25\\-0.5
    \end{bmatrix}, \quad V_2 \coloneqq \begin{bmatrix}
    -2\\5\\4\\3
    \end{bmatrix}.
\end{align*}
In particular, we have $\ell_v = \sigma^{\max}(V_{1})$ and in this special case, we can determine $\omega$ by inspection.

Hence, we can compute the bound on the Lipschitz constant of the constrained optimizer $\tilde{x}^\star$ as  $\ell_{t}=\sqrt{\frac{\beta}{\alpha}} \left(\frac{\ell_p}{\alpha} + \frac{\ell_C}{\omega}\right)$, and we can compute the analytical tracking bound as $\tfrac{\ell_t}{\alpha}\approx 0.5302$.

\begin{figure}[h!]
\centering
\includegraphics[width=.95\columnwidth]{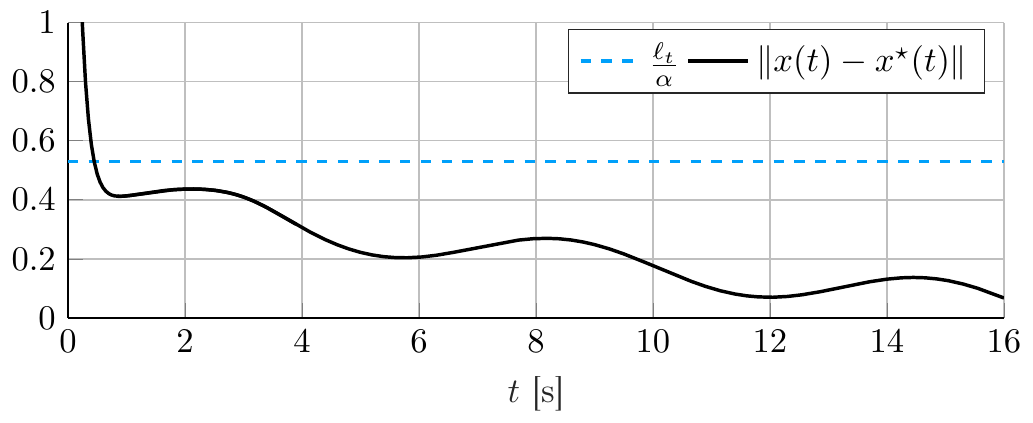}
\vspace{-.1em}
\caption{Tracking error $\norm{x(t)-{x}^\star(t)}$ (black, solid), and the tracking upper bound $\tfrac{\ell_t}{\alpha}$ (blue, dashed). \label{fig:sweeping.error}}
\end{figure}

\cref{fig:sweeping.error} presents the tracking error $\norm{x(t)-\txs(t)}$ (black, solid) and the tracking error upper bound $\frac{\ell_t}{\alpha}\approx 0.5302$ (blue, dashed). It can be seen that the tracking error is always within the bound, and that it is decreasing within the considered time interval.

\begin{figure}[tbh]
\vspace{-1em}
\includegraphics[width=0.95\columnwidth]{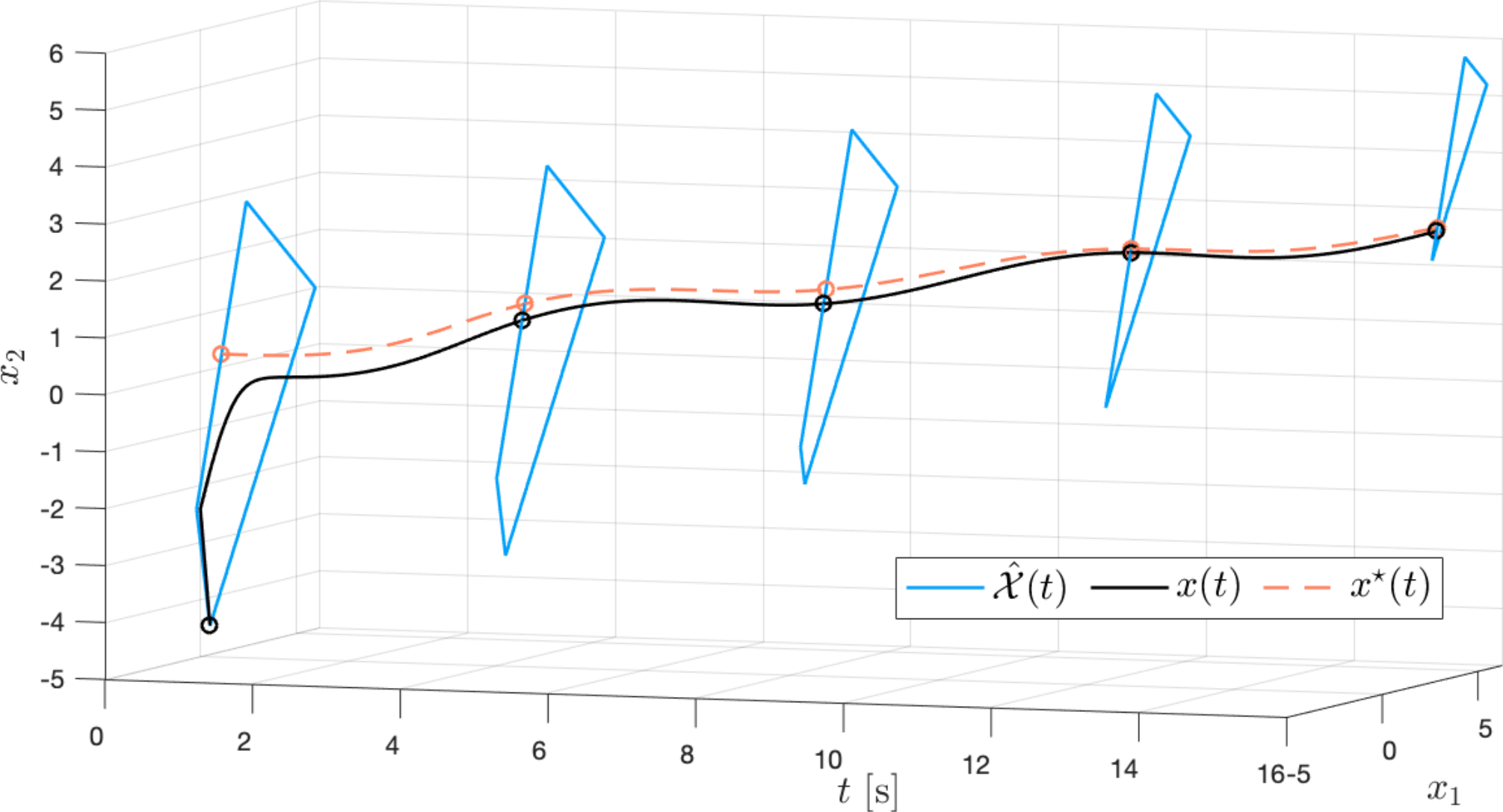}
\vspace{-.1em}
\caption{Two dimensional representation (i.e., $(x_1,x_2)$ plane) of the constrained time-varying optimization problem, captured in a different time steps, i.e, $t=0$s, $t=3.95$s, $t=7.95$s, $t=11.95$s, and $t=16$s. In each time step are presented: a constraint set $\hat{\calX} (t)$ (blue, solid), trajectory obtained via sweeping gradient flow, $x(t)$ (black, solid), instantaneous optimizer trajectory $x^\star(t)$ (vermilion, dashed). \label{fig:sweepingtrajectories}}
\end{figure}

In contrast to the unconstrained time-varying optimization, for the constrained time-varying optimization, the movement of both the cost function and the constraint set has to be considered. \cref{fig:sweepingtrajectories} illustrates the sweeping gradient flow trajectory $x(t)$ (black, solid), instantaneous optimizer trajectory $x^\star(t)$ (vermilion, dashed) and the time-varying constraint set $\hat{\calX} (t)$ (blue, solid), captured in a different time steps, i.e., $t=0$s, $t=3.95$s, $t=7.95$s, $t=11.95$s, and $t=16$s. Note, however, that contrary to the theoretical requirement, the problem becomes eventually infeasible with $\hat{\calX} (t)$ collapsing to an empty set.

\section{Conclusion}\label{sec:conc}
In this work, we have revisited classical sensitivity results to characterize and quantify the Lipschitz continuity  primal and dual solutions of nonlinear perturbed optimization problems. Since time-varying optimization is a natural application of our results, we have established a general tracking result for any continuous-time schemes satisfying a monotonicity-like property. In particular, this result can be immediately applied to discontinuous dynamical systems. To illustrate this fact, we have introduced \emph{sweeping gradient flows}, a novel discontinuous gradient descent scheme for convex, time-varying constrained optimization problems. This scheme is based on well-established perturbed sweeping processes and allows us to apply our generalized tracking result in a technically sophisticated context. Lastly, our numerical examples confirm our theoretical tracking guarantees.

It remains open whether some of our assumptions can be weakened: For instance, by considering the difference between cost function values, instead of the distance to the optimizer, it is presumably possible to establish tracking guarantees in terms of the objective value without requiring strong convexity. Similarly, there might be room to relax the LICQ assumption on the feasible set with more general constraint qualifications.

Finally, although our focus is on time-varying optimization, we envision our quantified sensitivity bounds to apply to other robustness questions about online optimization algorithms.

\appendix
\begin{lemma} \cite[Ch. 2.17]{bernsteinMatrixmathematicsTheory2009} \label{lem:block_inverse}
 Let $A \in \R^{n\times n}$, $B, C^\T \in \R^{n\times m}$ and $D \in \R^{m \times m}$. If $A$ is non-singular, then the  inverse
 \begin{align*}
  M^{-1}=\begin{bmatrix}
    		M_1 & M_2 \\ M_3 & M_4
    \end{bmatrix} \quad \text{of} \quad
    M\coloneqq \begin{bmatrix}
      A & C \\ B & D
    \end{bmatrix}
 \end{align*}
exists if and only if $D - BA^{-1}C$ is invertible. Then, we have
 \begin{align*}
      M_1  &= A^{-1} + A^{-1} C (D - BA^{-1} C)^{-1} B A^{-1}\,, \\
      M_2 &= -A^{-1} C (D - BA^{-1} C)^{-1} \,,\\
      M_3 &= - (D - BA^{-1} C)^{-1} B A^{-1}\,, \\
      M_4 & =  (D - BA^{-1} C)^{-1} \, .
  \end{align*}
\end{lemma}

\bibliographystyle{IEEEtran}
\bibliography{IEEEabrv,references}

\end{document}